\documentclass[twoside, 11pt]{article}
\usepackage{amssymb, amsmath, mathrsfs, amsthm}
\usepackage{graphicx}
\usepackage{color}
\usepackage{float, caption, subcaption}

\input{epsf}

\newcommand{\bd}{\begin{description}}
\newcommand{\ed}{\end{description}}
\newcommand{\bi}{\begin{itemize}}
\newcommand{\ei}{\end{itemize}}
\newcommand{\be}{\begin{enumerate}}
\newcommand{\ee}{\end{enumerate}}
\newcommand{\beq}{\begin{equation}}
\newcommand{\eeq}{\end{equation}}
\newcommand{\beqs}{\begin{eqnarray*}}
\newcommand{\eeqs}{\end{eqnarray*}}

\newcommand{\ceil}[1]{\left\lceil #1 \right\rceil}

\definecolor{DarkGreen}{rgb}{0.15, 0.4, 0.2}

\newcommand{\labelz}[1]{\label{#1}}

\newtheorem{theorem}{Theorem}

\newtheorem{lemma}{Lemma}

\newtheorem{case}{Case}
\newtheorem{subcase}{Subcase}[case]
\newtheorem{claim}{Claim}
\newtheorem{fact}{Fact}

\begin{document}
\title{\textbf{Ramsey and Gallai-Ramsey numbers for stars with extra independent edges} \footnote{Supported by the National Science Foundation of China (Nos. 11601254, 11551001, 11161037, 61763041, 11661068, and 11461054) and the Science Found of Qinghai Province (Nos.  2016-ZJ-948Q, and 2014-ZJ-907) and the  Qinghai Key Laboratory of Internet of Things Project (2017-ZJ-Y21).}
}

\author{
Yaping Mao\footnote{School of Mathematics and Statistis, Qinghai
Normal University, Xining, Qinghai 810008, China. {\tt
maoyaping@ymail.com}} \footnote{Academy of Plateau Science and Sustainability, Xining, Qinghai 810008, China}, \ \ Zhao Wang\footnote{College of Science,
China Jiliang University, Hangzhou 310018, China. {\tt
wangzhao@mail.bnu.edu.cn}}, \ \ Colton Magnant\footnote{Department
of Mathematics, Clayton State University, Morrow, GA, 30260, USA.
{\tt dr.colton.magnant@gmail.com}} \footnotemark[3], \ \ Ingo
Schiermeyer\footnote{Technische Universit{\"a}t Bergakademie
Freiberg, Institut f{\"u}r Diskrete Mathematik und Algebra, 09596
Freiberg, Germany. {\tt Ingo.Schiermeyer@tu-freiberg.de}}}

\maketitle

\begin{abstract}
Given a graph $G$ and a positive integer $k$, define the \emph{Gallai-Ramsey number} to be the minimum number of vertices $n$ such that any $k$-edge coloring of $K_n$ contains either a rainbow (all different colored) triangle or a monochromatic copy of $G$. In this paper, we obtain general upper and lower bounds on the Gallai-Ramsey numbers for the graph $G = S_t^{r}$ obtained from a star of order $t$ by adding $r$ extra independent edges between leaves of the star so there are $r$ triangles and $t - 2r - 1$ pendent edges in $S_t^{r}$. We also prove some sharp results when $t = 2$.
\end{abstract}

\section{Introduction}

In this work, we consider only edge-colorings of graphs. A coloring of a graph is called \emph{rainbow} if no two edges have the same color.

Edge colorings of complete graphs that contain no rainbow triangle have interesting and somewhat surprising structure. In 1967, Gallai \cite{MR0221974} examined this structure. His main result was restated in \cite{MR2063371} in the terminology of graphs and can also be traced back to \cite{MR1464337}. For the following statement, a trivial partition is a partition into only one part.

\begin{theorem}{\upshape \cite{MR1464337, MR0221974, MR2063371}}\labelz{Thm:G-Part}
In any coloring of a complete graph containing no rainbow triangle, there exists a nontrivial partition of the vertices (called a Gallai-partition) such that there are at most two colors on the edges between the parts and only one color on the edges between each pair of parts.
\end{theorem}

The induced subgraph of a Gallai colored complete graph constructed by selecting a single vertex from each part of a Gallai partition is called the \emph{reduced graph}. By Theorem~\ref{Thm:G-Part}, the reduced graph is a $2$-colored complete graph.

Given two graphs $G$ and $H$, let $R(G, H)$ denote the $2$-color Ramsey number for finding a monochromatic $G$ or $H$, that is, the minimum number of vertices $n$ needed so that every red-blue coloring of $K_{n}$ contains either a red copy of $G$ or a blue copy of $H$. Although the reduced graph of a Gallai partition uses only two colors, the original Gallai colored complete graph could certainly use more colors. With this in mind, we consider the following generalization of the Ramsey numbers. Given two graphs $G$ and $H$, the \emph{general $k$-colored Gallai-Ramsey number} $gr_k(G:H)$ is defined to be the minimum integer $m$ such that every $k$-coloring of the complete graph on $m$ vertices contains either a rainbow copy of $G$ or a monochromatic copy of $H$. Since we have the additional restriction of forbidding the rainbow copy of $G$, it is clear that $gr_k(G:H)\leq R_k(H)$ for any graph $G$. Gallai-Ramsey numbers have been studied for a wide variety of monochromatic graphs. We refer the intereted reader to \cite{MR2606615} for a survey of relevant results with an updated version available at \cite{FMO14}.

The graph $S_t^{r}$ is obtained from a star of order $t$ by adding an extra $r$ independent edges between the leaves of the so that there are $r$ triangles and $t - 2r - 1$ pendent edges in $S_t^{r}$. For $r=0$ we obtain $S_t^{r} = K_{1,t-1}$, which are called \emph{stars}. For $r = \frac{t-1}{2}$, if $t$ is odd we obtain $S_t^{r} = F_{\frac{t-1}{2}}$, which are called \emph{fans}. Gallai-Ramsey numbers for stars and fans have been considered in \cite{MR2063371} and \cite{Fans}, respectively.

Therefore, in this paper we deal with those graphs $S_t^{r}$ where $1 \leq r < \frac{t-1}{2},$ i.e. where $S_t^{r}$ is neither a star nor a fan. We prove general bounds on the Gallai-Ramsey number for $S_{t}^{r}$ and sharp results for some cases when $t = 2$.

Our first main result finds the $2$-color Ramsey numbers for $S_{t}^{r}$ for a variety of situations.

\begin{theorem}\labelz{Thm:St2Classic}
$(1)$ For $t\geq 7$, $R(S_t^{2}, S_t^{2})=2t-1$.

$(2)$ For $t\geq 15$, $R(S_t^{3}, S_t^{3})=2t-1$.

$(3)$ For $t\geq 6r-5$, $R(S_t^{r}, S_t^{r})=2t+2r-1$.
\end{theorem}

For the graphs $S_t^{2}$, we have the following.

\begin{theorem}\labelz{Thm:St2}
$(1)$ For $k\geq 1$,
$$
gr_k(K_3;S_6^{2})=\begin{cases}
2\times 5^{\frac{k}{2}}+\frac{1}{4}\times 5^{\frac{k-2}{2}}+\frac{3}{4}, &\mbox {\rm if}~k~is~even;\\[0.2cm]
\lceil\frac{51}{10}\times 5^{\frac{k-1}{2}}+\frac{1}{2}\rceil,
&\mbox {\rm if}~k~is~odd.
\end{cases}
$$

$(2)$ For $k\geq 3$,
$$
gr_k(K_3;S_8^{2})=\begin{cases}
14\times 5^{\frac{k-2}{2}}+\frac{1}{2}\times 5^{\frac{k-4}{2}}+\frac{1}{2}, &\mbox {\rm if}~k~is~even;\\[0.2cm]
7\times 5^{\frac{k-1}{2}}+\frac{1}{4}\times
5^{\frac{k-3}{2}}+\frac{3}{4}, &\mbox {\rm if}~k~is~odd.
\end{cases}
$$

$(3)$ For $k\geq 1$ and $t\geq 6$,
$$
\begin{cases}
2(t-1)\times 5^{\frac{k-2}{2}}+1 \leq gr_k(K_3;S_t^{2}) \leq 2t\times 5^{\frac{k-2}{2}}, &\mbox {\rm if}~k~is~even;\\[0.2cm]
(t-1)\times 5^{\frac{k-1}{2}}+1 \leq gr_k(K_3;S_t^{2}) \leq t\times
5^{\frac{k-1}{2}}, &\mbox {\rm if}~k~is~odd.
\end{cases}
$$
\end{theorem}

Finally, we also provide general bounds on the Gallai-Ramsey numbers for $S_{t}^{r}$.

\begin{theorem}\labelz{Thm:Str}
For $t\geq 6r-5$,
$$
\begin{cases}
2(t-1)\times 5^{\frac{k-2}{2}}+1\leq gr_k(K_3;S_t^{r})\leq [2t+8(r-1)]\times 5^{\frac{k-2}{2}}-4(r-1), \\
\hspace{3.5in} \text{ if $k$ is even;}\\
(t-1)\times 5^{\frac{k-1}{2}}+1\leq gr_k(K_3;S_t^{r})\leq [t+4(r-1)]\times 5^{\frac{k-1}{2}}-4(r-1), \\
\hspace{3.5in} \text{ if $k$ is odd.}
\end{cases}
$$
\end{theorem}

In Section~\ref{Sec:Classical}, we prove Theorem~\ref{Thm:St2Classic}. The proof of Theorem~\ref{Thm:St2} is mostly contained in Section~\ref{Sec:GR}. We omit the proofs of Item~$(3)$ of Theorem~\ref{Thm:St2} and Theorem~\ref{Thm:Str} since they follow a similar structure to the proofs of the other parts of Theorem~\ref{Thm:St2}. The complete proofs will be made publicly available.

\section{Results for the classical Ramsey number}\label{Sec:Classical}

The following lemma is almost immediate.

\begin{lemma}\labelz{Lemma:MonoSmallParts}
If $G$ is a Gallai colored complete graph of order at least $4n-3$ in which all parts of a Gallai partition have order at most $n - 1$ and all edges in between the parts of $G$ have one color, say red, then $G$ contains a red copy of $F_{n}$.
\end{lemma}

\begin{proof}
Let $H_{1}, H_{2}, \dots, H_{t}$ be the parts of the assumed partition, so since $|G| \geq 4n - 3$, we see that $t \geq 5$. Since $|H_{i}| \leq n - 1$, there exists an integer $r$ and corresponding set of parts $H_{2}, H_{3}, \dots, H_{r}$ such that $n \leq |H_{2} \cup H_{3} \cup \dots \cup H_{r}| \leq 2n - 2$. This, in turn, implies that $|H_{r + 1} \cup H_{r + 2} \cup \dots \cup H_{t}| \geq n$. Then a single vertex from $H_{1}$ along with $n$ red edges from $H_{2} \cup H_{3} \cup \dots \cup H_{r}$ to $H_{r + 1} \cup H_{r + 2} \cup \dots \cup H_{t}$ produces a red copy of $F_{n}$.
\end{proof}

\subsection{The case $r=2$}

For $r=2$, we need to prove that $R(S_t^{2}, S_t^{2})=2t-1$.

~

\noindent {\it Proof of $(1)$ of Theorem \ref{Thm:St2Classic}}. For the lower bound, let $G$ be a $2$-edge colored graph obtained from two copies of $K_{t-1}$ by adding all blue edges between them. Clearly, there is neither a red $S_t^{2}$ nor a blue $S_t^{2}$ in $G$. So $R(S_t^{2}, S_t^{2})\geq 2t-1$.

For the upper bound, let $G$ be a $2$-colored copy of $K_{2t - 1}$ and for a contradiction, suppose $G$ contains no monochromatic copy of $S_{t}^{2}$. For each $v\in V(G)$, let $A_v$ and $B_v$ be the set of vertices with red and blue edges respectively to $v$. Suppose that there exists a vertex $v\in V(G)$ with $|A_v|\geq t+1$.

\begin{fact}\labelz{Fact:NoMatchingInAorB} 
There are neither a red copy of $2K_{2}$ within $A_v$ nor a blue copy of $2K_{2}$ within $B_v$.
\end{fact}

By Fact~\ref{Fact:NoMatchingInAorB}, there is at most a red triangle or a red star in $A_v$. If there is a red star with center $u$ in $A_v$, then $A_v-u$ contains a blue $K_t$, and hence $A_v-u$ contains a blue $S_t^{2}$, a contradiction. If there is a red triangle in $A_v$, then $A_v$ is a blue graph obtained from a $K_{t+1}$ by deleting a triangle. clearly, $A_v$ contains a blue $S_t^{2}$, a contradiction.

Suppose that there exists a vertex $v\in V(G)$ such that $|A_v|=t$. Then $|B_v|=t-2$. Similar to Fact~\ref{Fact:NoMatchingInAorB}, there is no red copy of $2K_{2}$ within $A_v$, and hence there is at most a red triangle or a red star in $A_v$. If there is a red triangle in $A_v$, then $A_v$ contains a blue graph obtained from $K_{t-3}$ and a star $K_{1,3}$ by identifying the center of the star and a vertex of $K_{t-3}$. Clearly, there is a blue $S_t^{2}$, a contradiction. If there is a red star with center $u$ in $A_v$, then we let $A_v'=A_v-u$. Clearly, $A'=K_{t-1}$, and hence each edge from $u$ to $A_v'$ must be red. To avoid a blue $S_t^{2}$, the edges from $A_v'$ to $B$ must be red. If there is a red edge in $B$, then the edges among $v,u,a,B$ form a red $S_t^{2}$, where $a\in A_v'$, a contradiction. So $B$ is a blue clique of order $t-2$. If there is an red edge from $u$ to $B$, then $u,v,A'$ form a red $S_t^{2}$, a contradiction. Therefore, the edges from $u$ to $B$ are red, and hence $u,v,B$ form a blue $S_t^{2}$, a contradiction.

For the remainder of the proof, we need only the following fact.

\begin{fact}\labelz{Fact:t-1} 
For each $v\in V(G)$, $|A_v|=|B_v|=t-1$.
\end{fact}

From Fact~\ref{Fact:NoMatchingInAorB}, there is neither a red copy of $2K_{2}$ within $A_v$ nor a blue copy of $2K_{2}$ within $B_v$. Therefore, there is at most a red triangle or a red star in $A_v$, and there is at most a blue triangle or a blue star in $B_v$.

Suppose that there is a red star with center $a$ in $A_v$ and there is a blue star with center $b$ in $B_v$. Let $A_v'=A_v-a$ and $B_v'=B_v-b$. Then $A_v'$ is a red clique of order $t-2$, and $B_v'$ is a blue clique of order $t-2$. If $t\geq 7$, then to avoid a red $S_t^{2}$, the red edges from $A'$ to $B'$ form a red matching. Similarly, the blue edges from $A'$ to $B'$ form a blue matching. Since $t\geq 7$, it follows that the number of red and blue edges is $(t-2)^2>2(t-2)$, a contradiction.

Suppose that there is a red triangle $x_{A}y_{A}z_{A}x_{A}$ in $A_v$ and there is a blue triangle $x_{B}y_{B}z_{B}x_{B}$ in $B_v$. Let $A_v'=A_v-\{x_{A},y_{A},z_{A}\}$ and $B_v'=B_v-\{x_{B},y_{B},z_{B}\}$. It is clear that $A_v'$ is a blue clique of order $t-4$ and $B_v'$ is a red clique of order $t-4$. Choose a vertex $w\in A_v'$. From Fact 2, there are $t-2$ red edges and one blue edge from $w$ to $B_v$ since $|B|=t-1$. If $t\geq 7$, then there is a red $2$-matching in $B$, and hence there is a red $S_t^{2}$, a contradiction.

Suppose that there is a red star with center $a$ in $A_v$ and there is a blue triangle $xyzx$ in $B_v$. Let $A_v'=A_v-a$ and $B_v'=B_v-\{x,y,z\}$. Clearly, $A_v'$ is a blue clique of order $t-2$ and $B_v'$ is a red clique of order $t-4$. If there is a blue edge from $A'$ to $B$, then this edge together with $A_v'$ form a blue $S_t^{2}$, a contradiction. So we can assume that the edges from $A_v'$ to $B_v$ are red. Since $t\geq 7$, it follows that the edges from $A_v'$ to $B_v$ and $B$ from a red $S_t^{2}$, a contradiction. \qed

\subsection{The case $r=3$}

For $r=3$, we need to prove that $R(S_t^{3}, S_t^{3})=2t-1$. 

~

\noindent {\it Proof of $(2)$ of Theorem \ref{Thm:St2Classic}}. The lower bound follows from the same example as presented in the proof of $(1)$.  For the upper bound, let $G$ be a $2$-coloring of $K_{2t - 1}$ and for a contradiction, suppose $G$ contains no monochromatic copy of $S_{t}^{3}$.  For each $v\in V(G)$, let $A_v$ and $B_v$ be the set of vertices incident red and blue edges to $v$, respectively.

\setcounter{case}{0}
\begin{case}
There exists a vertex $v\in V(G)$ such that $|A_v|\geq t+2$.
\end{case}

It is clear that $A$ contains no red $3$-matching (that is $3K_{2}$). Let $X=\{u_1,u_2,u_3,u_4\}$ be the vertices of a red $2$-matching, and let $A'=A-X$. Clearly, $A'$ is a blue clique of order $t-2$. To avoid a red $3$-matching, for each vertex $w$ in $A'$, the edges from $w$ to $X$ are red or there is at most one blue edge from $w$ to $X$. Choose $Y=\{w_1,w_2,w_3,w_4\}$. Then the edges between $X$ and $Y$ form a red subgraph obtained from $K_{4,4}$ by deleing at most four edges, and hence there are a $3$-matching in this red subgraph, a contradiction.

\begin{case}
There exists a vertex $v\in V(G)$ such that $|A_v|=t+1$.
\end{case}

To avoid a red $S_t^3$, $A$ contains at most a red $2$-matching, say $\{u_1u_2,u_3u_4\}$. Let $A'=A-\{u_1,u_2,u_3,u_4\}$.

\begin{claim}\labelz{Claim:ui}
For each $u_i \ (1\leq i\leq 4)$, the edges from $u_i$ to $A'$ are red or blue.
\end{claim}

\begin{proof}
Assume, to the contrary, that there exist two vertices $w_1,v_1$ in $A'$ such that $u_1w_1$ is blue and $u_1v_1$ is red. To avoid a red $3$-matching, $u_2w_1$ is blue and for $x\in A'-\{w_1,v_1\}$, $u_2x$ is blue. To avoid a blue $S_{t}^3$, $u_3w_1$ and $u_4w_1$ are red. Clearly, $xu_3$ and $xu_4$ are blue. Then there is a blue $S_{t}^3$, a contradiction.
\end{proof}

Suppose that the edges from $\{u_1,u_2\}$ to $A'$ are red. Since $t\geq 15$, there exist two vertices $w_1,w_2\in A'$ such that $\{u_1w_1,u_2w_2\}$ is a red $2$-matching, and hence $\{u_1w_1,u_2w_2,u_3u_4\}$ is a red $3$-matching, a contradiction. We can also get a contradiction if the edges from $\{u_3,u_4\}$ to $A'$ are red.

Suppose that the edges from $\{u_1,u_2\}$ to $A'$ are blue. To avoid a blue $S_{t}^3$, the edges from $\{u_3,u_4\}$ to $A'$ are red, and hence there is a red $3$-matching in $A$, a contradiction. We can also get a contradiction if the edges from $\{u_3,u_4\}$ to $A'$ are blue.

From now on, without loss of generality, we assume that the edges from $\{u_1,u_3\}$ to $A'$ are red and the edges from $\{u_2,u_4\}$ to $A'$ are blue. To avoid a red $3$-matching in $A=A'\cup \{u_1,u_2,u_3,u_4\}$, $u_2u_4$ is blue. Let $A''=A'\cup \{u_2,u_4\}$. It is clear that the graph induced by $A''$ is a blue clique of order $t-1$. To avoid a blue $S_t^3$, each edge from $A''$ to $B$ is red. The following facts are immediate.

\begin{fact}\labelz{Fact:no3}
There is no vertex incident to at least $3$ red edges in $B$, and there is no red $2$-matching in $B$.
\end{fact}

From Fact~\ref{Fact:no3}, there is at most one red triangle in $B$.

\begin{claim}\labelz{Claim:oneRed}
There is at most one red edge from $u_1$ (resp. $u_2$) to $B$.
\end{claim}

\begin{proof}
Assume, to the contrary, that there are two edges $u_1w_1,u_1w_2$ from $u_1$ to $\{w_1,w_2\}$, where $w_1,w_2\in B$. Choose $u\in A''-\{u_1,u_4\}$. Then $uvu_1v$, $u_1u_2w_1u_1$ and $u_1u_4w_2u_1$ form $3$ red triangles and hence there is a red $S_{t}^3$, a contradiction.
\end{proof}

Similarly, there is at most one red edge from $u_2$ to $B$. Let $h\in B-\{w_1,w_2,x,y,z\}$. Clearly, $hu_1,hu_2,hv$ are blue, and hence there is a blue $S_{t}^3$, a contradiction.

\begin{case}
For each vertex $v\in V(G)$, $|A_v|=|B_v|=t-1$.
\end{case}

Since $t\geq 15$ and $R(F_3,F_3)=13$, it follows that there is a red $F_3$ or a blue $F_3$ in $G$. Without loss of generality, let $F_3$ be a red fan of order $7$, and let $v$ be the center of $F_3$. In this case, there $t-7$ pendent edges incident to $v$. There is a red $S_{t}^3$, a contradiction.

\begin{case}
There exists a vertex $v\in V(G)$ such that $|A_v|=t$.
\end{case}

It is clear that there is at most a red $2$-matching in $A$. Let $X=\{u_1,u_2,u_3,u_4\}$, and let $A'=A-X$. Then $|A'|=t-4$ and $A'$ be a blue clique of order $t-4$.

\begin{fact}\labelz{Fact:3blue}
For each $w\in A'$, there are at most $3$ blue edges from $w$ to $X$.
\end{fact}

We distinguish the following cases to show our proof.

\begin{subcase}
There exists a vertex $w\in A'$ such that there are exactly three blue edges from $w$ to $X$.
\end{subcase}

To avoid a blue $S_t^3$, the edges from $w$ to $X$ are red. To avoid a red $3$-matching, for any $x\in A'-w$, the edges from $w$ to $X$ are blue, and hence there is a blue $S_t^3$, a contradiction.

\begin{subcase}
There exists a vertex $w\in A'$ such that there are exactly two blue edges from $w$ to $X$.
\end{subcase}

Let $wv_1$ and $wv_2$ be two blue edges from $w$ to $X$, and let $B'=B-\{v_1,v_2\}$. To avoid a blue $S_t^3$, there are three edges, say $wu_1,wu_2,wu_3$, from $w$ to $X$. To avoid a red $3$-matching, the edges from $x$ to $X-v_3$ are blue.

Suppose that there are at least two red edges $u_1v',u_1v''$ from $u_1$ to $B$. Then we have the following claim. 

\begin{claim}\labelz{Claim:u2blue}
The edges from $u_2$ to $B$ are blue.
\end{claim}

\begin{proof}
Assume, to the contrary, that there exists a vertex $v'''$ (note that $v'''$ is not necessarily different from $v'$ and $v''$). Then there is a red $S_t^3$, a contradiction. Then there is a red $S_t^3$, where $uwu_1v,wu_1v'w,wu_2v'''w$ are three triangles in $S_t^3$, a contradiction.
\end{proof}

From Claim~\ref{Claim:u2blue}, the graphs induced by the vertices in $A'\cup B\cup \{u_2\}$ contains a blue $S_t^3$, a contradiction.

\begin{subcase}
For each $w\in A'$, the edges from $w$ to $B$ are red.
\end{subcase}

Note that $|A'|=t-4$ and $|B|=t-2$.

\begin{fact}\labelz{Fact:3red}
There is no vertex incident at least $3$ red edges, and there is no
red $2$-matching in $B$.
\end{fact}

From Fact~\ref{Fact:3red}, there is at most a red triangle $xyzx$ in $B$. Let $B'=B-\{x,y,z\}$. We claim that the edges from $X$ to $B'$ are red. Assume, to the contrary, that there exists a vertex $w\in B'$ such that $wu_1$ is blue. Since $vw,wx,wy,wz$ are blue, it follows that there is a blue $S_t^3$, a contradiction.

Choose $w_1,w_2\in A'$ and $v_1,v_2\in B'$. Then there exists a red $S_t^3$, $u_1wv_1u_1,u_1w_2v_2u_1,u_1u_2vu_1$ are three triangles in $S_t^3$, a contradiction. \qed

\subsection{For general $r$}

For general $r$, we need to prove that $R(S_t^{r}, S_t^{r})=2t-1$ for $t\geq 6r-5$.

~

\noindent {\it Proof of $(3)$ of Theorem~\ref{Thm:St2Classic}}. The lower bound follows from the same example as presented in the proof of $(1)$.  For the upper bound, let $G$ be a $2$-coloring of $K_{2t - 1}$ and for a contradiction, suppose $G$ contains no monochromatic copy of $S_{t}^{3}$. For each $v\in V(G)$, let $A_v$ and $B_v$ be the set of vertices incident red and blue edges to $v$, respectively.

\setcounter{case}{0}
\begin{case}
There exists a vertex $u\in V(G)$ such that $|A_u|\geq t$ and there exists a vertex $v\in V(G)$ such that $|A_v|\geq t$.
\end{case}

To avoid a red $S_t^r$, there are at most $r-1$ red matching in $A_u$. By deleting $2r-2$ vertices in $A_u$, the resulting graph has no red edges and contains a blue clique $A'=K_{t-2r+2}$. Similarly, by deleting $2r-2$ vertices in $A_v$, the resulting graph has no blue edges and contains a red clique $B'=K_{t-2r+2}$. Choose $x\in A'$. To avoid a blue $S_t^r$, the number of blue edges from $x$ to $B'$ is at most $2r-3$, and hence the number of red edges from $x$ to $B'$ is at least $t-4r+5$. Since $t\geq 6r-5$, it follows that there is a red fan $F_{2r}$ with center $x$, and there is a blue fan $F_{2r}$ with center $x$. Since the degree of $x$ is $2t-2$, it follows that the red or blue degree of $x$ is $t-1+(2r+1)=t+2r$, a contradiction.

\begin{case}
For each vertex $v\in V(G)$, $|A_v|\geq t$ or $|B_v|\geq t$.
\end{case}

Without loss of generality, we suppose $|A_v|\geq t$. Then there is a red Hamilton cycle $C$ in $G$. If there exists a vertex $x$ of degree at least $t+r$, then there is a red fan $F_{2r}$, and hence there is a red $S_t^r$, a contradiction. Suppose that for each vertex $x$ its red degree is at most $t+r-1$ and at least $t$.

\begin{claim}
$(3')$ For $t\geq 6r-4$, $R(S_t^{r}, S_t^{r}) \leq 2t+2r-1$.
\end{claim}

\begin{proof}
Let $G$ be a graph of order $2t+2r-1,$ whose edges are colored red and blue. For each $v\in V(G)$, let $A_v$ and $B_v$ be the set of vertices incident with red and blue edges to $v$, respectively. We distinguish the following two cases, which cover all possible colorings.

\setcounter{case}{0}
\begin{case}
There exists a vertex $u\in V(G)$ such that $|A_u|\geq t-1$ and there exists a vertex $v\in V(G)$ such that $|B_v|\geq t-1$.
\end{case}

We follow the proof of Case $1$ above with some small changes.

To avoid a red $S_t^r$, there is at most a $r-1$ red matching in $A_u$. By deleting $2r-2$ vertices in $A_u$, the resulting graph has no red edges and contains a blue clique $A'=K_{t-2r+1}$. Similarly, by deleting $2r-2$ vertices in $A_v$, the resulting graph has no blue edges and contains a red clique $B'=K_{t-2r+1}$. Choose $x\in A'$. To avoid a blue $S_t^r$, the number of blue edges from $x$ to $B'$ is at most $2r-2$, and hence the number of red edges from $x$ to $B'$ is at least $t-4r+4$. Since $t\geq 6r-4$, it follows that there is a red fan $F_{2r}$ with center $x$, and there is a blue fan $F_{2r}$ with center $x$. Since the degree of $x$ is $2t+2r-2$, it follows that the red or blue degree of $x$ is at least $t+r-1$, a contradiction.

\begin{case}
For each vertex $v\in V(G)$, $|A_v|\geq t+2r$ or $|B_v|\geq t+2r$.
\end{case}

Without loss of generality, we suppose $|A_v|\geq t+2r$. Let $\delta^R(G)$ denote the red minimum degree in $G.$ Then there is a red Hamilton cycle $C$ in $G$, since $\delta^R(G) \geq t+2r > \frac{2t+2r-1}{2}$ (Dirac). Since  $|A_v|\geq t+2r$ for every vertex $v \in V(G),$ we have $d^R(v) > \frac{2t+2r-1}{2} + r.$ Using a combinatorial counting argument we conclude that every vertex $v \in V(G)$ is contained in a red fan $F_{2r}$, and hence there is a red $S_t^r$, a contradiction.
\end{proof}

This completes the proof of Theorem~\ref{Thm:St2Classic}. \qed

\section{Results for Gallai Ramsey number}\label{Sec:GR}

In this section, we study the Gallai Ramsey number of stars with extra independent edges.


\subsection{The case $r=2$ and $t=6$}

We first give a lower bound construction for item $(1)$ in Theorem~\ref{Thm:St2}.

\begin{lemma}\labelz{Lem:Lower-S62}
For $k \geq 1$, we have
$$
gr_k(K_3;S_6^{2})\geq \begin{cases}
2\cdot 5^{\frac{k}{2}}+\frac{1}{4}\cdot 5^{\frac{k-2}{2}}+\frac{3}{4}, & \text{ if $k$ is even;}\\
\ceil{ \frac{51}{10}\cdot 5^{\frac{k-1}{2}}+\frac{1}{2} }, & \text{ if $k$ is odd.}
\end{cases}
$$
\end{lemma}

\begin{proof}
Let $G_1$ be a copy of $K_5$ entirely colored by color $1$ and let $G_2$ be a copy of $K_{10}$ consisting of two copies of $G_1$ joined by all edges of color $2$. Now suppose we have constructed a Gallai colored complete graph $G_{k-2}$ with $k \geq 3$ using $k-2$ colors which contains no monochromatic copy of $S_6^2$. We construct the graph $G_k$ by making five copies of $G_{k-2}$ and inserting edges of colors $k-1$ and $k$ between the copies to form a blow-up of the unique $2$-coloring of $K_5$ with no monochromatic triangle. This coloring contains no rainbow triangle and no monochromatic copy of $S_6^2$, and has order
$$
|G_{k}| = \begin{cases}
2 \cdot 5^{k/2} & \text{ if $k$ is even;}\\
5^{(k+1)/2} & \text{ if $k$ is odd.}
\end{cases}
$$
For $k=3$ we set $G'_3 = G_3$.

For all $k \geq 4$ we extend this construction as follows. For $j \geq 4$ define the graph $A_1^j$ to be a colored complete graph $K_6$ containing a perfect matching using colors $2, 3$ and $j$, with all other edges in color $1$. If $k$ is odd, then choose two copies of $G'_{k-2}$ and replace one copy of $G_1$ within each of them by a copy of $A_1^{k-1}$ and $A_1^k$, respectively. This modification process is applied at each step of the construction above so in particular,
\bi
\item $|G'_{3}| = 25$,
\item $|G'_{5}| = 5 \cdot 25 + 2 = 127$,
\item $|G'_{7}| = 5 \cdot 127 + 2 = 637$,
\ei
and so on. This yields
$$
|G_k| = 5\left(\ceil{ \frac{51}{10} \cdot 5^{\frac{k-3}{2}}+\frac{1}{2} } - 1\right)+2 =\ceil{ \frac{51}{10} \cdot 5^{\frac{k-1}{2}}-\frac{1}{2} },
$$
as claimed.

For $j \geq 4$ let $B_1^j$ be a complete graph $K_6$ containing a perfect matching using colors $1, j-1$ and $j$, and all other edges in color $1$. Now if $k$ is even, then choose one copy of $G'_{k-2}$ and replace one copy of $G_1$ in it by a copy of $B_1^k$. This modification process is applied at each step of the construction above so in particular,
\bi
\item $|G'_{2}| = 10$,
\item $|G'_{4}| = 5 \cdot 10 + 1 = 51$,
\item $|G'_{6}| = 5 \cdot 51 + 1 = 256$,
\ei
and so on. This yields
$$
|G_k| = 5\left(\frac{41}{4} \cdot 5^{\frac{k-4}{2}}+\frac{3}{4} - 1\right)+1 = \frac{41}{4} \cdot 5^{\frac{k-2}{2}}-\frac{1}{4},
$$
as claimed.
\end{proof}

Next, we give the upper bound for $(1)$ of Theorem~\ref{Thm:St2}.

\begin{lemma}\labelz{Lem:Upper-S62}
$$
gr_k(K_3;S_6^{2})\leq \begin{cases}
2\cdot 5^{\frac{k}{2}}+\frac{1}{4}\cdot 5^{\frac{k-2}{2}}+\frac{3}{4}, & \text{ if $k$ is even;}\\
\ceil{ \frac{51}{10}\cdot 5^{\frac{k-1}{2}}+\frac{1}{2} }, & \text{ if $k$ is odd.}
\end{cases}
$$
\end{lemma}

\begin{proof}
We prove this upper bound by induction on $k$. The case $k=1$ is trivial and the case $k=2$ is precisely the statement $R(S_6^{2},S_6^{2})=15$. We therefore suppose $k\geq 3$ and let $G$ be a coloring of $K_n$ containing no rainbow triangle and no monochromatic copy of $S_6^{2}$ where
$$
n= n(k) = \begin{cases}
2\cdot 5^{\frac{k}{2}}+\frac{1}{4}\cdot 5^{\frac{k-2}{2}}+\frac{3}{4}, & \text{ if $k$ is even;}\\
\ceil{ \frac{51}{10}\cdot 5^{\frac{k-1}{2}}+\frac{1}{2} }, & \text{ if $k$ is odd.}
\end{cases}
$$

Since $G$ is a Gallai coloring, it follows from Theorem~\ref{Thm:G-Part} that there is a Gallai partition of $V(G)$. Suppose that the two colors appearing in the Gallai partition are red and blue. Let $m$ be the number of parts in this partition and choose such a partition where $m$ is minimized. Let $H_{1}, H_{2}, \dots, H_{m}$ be the parts of this partition, say with $|H_{1}| \geq |H_{2}| \geq \dots \geq |H_{m}|$. Let $q$ be the number of ``large'' parts of order at least $3$. This means that $|H_{q}| \geq 3$ and $|H_{q + 1}| \leq 2$. Let $X$ be the set of vertices in the ``small'' parts of order at most $2$. Note that $m \leq 10$ by Theorem~\ref{Thm:St2Classic}. First some easy observations related to facts noted in \cite{Fans}. The proof of each involves constructing a monochromatic copy of $S_{6}^{2}$ in the assumed structure.

\begin{fact}\labelz{Fact:Structures}
Let $A$ and $B$ be two disjoint non-empty sets of vertices in a colored complete graph with no monochromatic copy of $S_{6}^{2}$, with all edges from $A$ to $B$ colored red.
\bi
\item If $|B| \geq 3$, then $A$ contains no vertex with two incident red edges (within $A$).
\item If $|A| \geq 5$, then $A$ contains no red copy of $2K_{2}$.
\item If $|A| \geq 4$, then at most one of $A$ or $B$ contains a red edge.
\ei
Furthermore, given three disjoint non-empty sets of vertices with all red edges between pairs of these sets, then either all three of these sets must have order at most $2$ or two of these sets must have order $1$.
\end{fact}

We break the proof into two cases, where $k = 3$ and where $k \geq 4$.

\setcounter{case}{0}
\begin{case}
$k = 3$.
\end{case}

With $k = 3$, we have $n=26$ and say the three colors used are red, blue, and green.

If $2\leq m\leq 3$, then by the minimality of $m$, we may assume $m=2$, say with corresponding parts $H_1$ and $H_2$ with all red edges in between the two parts. Since $n=26$, it follows that $|H_{1}|\geq 13$. By Fact~\ref{Fact:Structures}, $H_1$ does not contain red $2K_2$. Thus, by deleting at most $2$ vertices from $H_{1}$, we can obtain a subgraph of $H_{1}$ in which there are no red edges. Since $R(S_6^2,S_6^2)=11$, it follows that there is a blue or green copy of $S_6^2$ within this subgraph, a contradiction. We may therefore assume that $4 \leq m \leq 10$.

If $q \geq 6$, then there is a monochromatic triangle in the reduced graph restricted to the $q$ large parts, contradicting Fact~\ref{Fact:Structures}. At the opposite extreme, if $q = 0$, then all parts have order at most $2$ so there are at least $13$ parts. Since $R(S_6^2,S_6^2)=11$, it follows that there is a monochromatic copy of $S_6^2$ in the reduced graph, a contradiction. Thus, we may assume that $1 \leq q \leq 5$.

If $q=5$, then by Fact~\ref{Fact:Structures}, in order to avoid a monochromatic triangle in the reduced graph that uses at least two vertices corresponding to large parts, we must have $m = q = 5$ and the parts are arranged to form a blow-up of the unique $2$-coloring of $K_{5}$ with no monochromatic triangle. This means that $|H_{1}|\geq \ceil{ \frac{26}{5} } =6$. By Fact~\ref{Fact:Structures}, $H_{1}$ contains no red or blue copy of $2K_{2}$ and no vertex with two incident red edges or two incident blue edges. This means that $H_{1}$ contains at most one red edge and at most one blue edge, meaning that $H_{1}$ contains a green copy of $S_{6}^{2}$, for a contradiction. More generally, this argument using Fact~\ref{Fact:Structures} yields the following fact.

\begin{fact}\labelz{Fact:AtMost5}
With $k = 3$, if a part $A$ of the Gallai partition has all red edges to at least $3$ vertices and all blue edges to at least $3$ vertices, then $|A| \leq 5$.
\end{fact}

Next suppose $q = 4$. In order to avoid a monochromatic triangle among the large parts, there are only two possible colorings for the reduced graph restricted to the large parts, one with complementary monochromatic copies of $P_{4}$ and the other with a $C_{4}$ in one color and $2K_{2}$ in the other. In either case, each large part has red and blue edges to another large part so by Fact~\ref{Fact:AtMost5}, there are at most $20$ vertices in large parts so $|X| \geq 6$. If the reduced graph restricted to the large parts contains a $C_{4}$ in one color, then the edges from $X$ to the large parts cannot be colored without creating a monochromatic copy of $S_{6}^{2}$, so we may assume the reduced graph restricted to the large parts contains two complementary monochromatic copies of $P_{4}$. In order to avoid a monochromatic triangle in the reduced graph using at least two large parts, the edges from all of $X$ to the ends of each $P_{4}$ must have the same color as the $P_{4}$ itself. By the minimality of $m$, $X$ must be a single part of the Gallai partition of order at least $6$, a contradiction.


Next suppose $q=3$. By Fact~\ref{Fact:Structures}, we may assume that the edges from $H_1$ to $H_2\cup H_3$ are red and the edges from $H_2$ to $H_3$ are blue. If a part in $X$ had red edges to $H_{1}$, then to avoid a monochromatic triangle in the reduced graph, it must have blue edges to both $H_{2}$ and $H_{3}$, creating a monochromatic copy of $S_{6}^{2}$. Thus, the edges from $X$ to $H_1$ must all be blue. By Fact~\ref{Fact:Structures}, $|H_{2}|\leq 5$ and $|H_{3}|\leq 5$. Then we have the following claim.

\begin{claim}\labelz{Clm:5}
There is no set of at least $7$ vertices in $X$ with all one colored edges to a large part.
\end{claim}

In particular, since $H_{1}$ has all blue edges to $X$, this claim implies that $|X|\leq 6$.

\begin{proof}
Assume, to the contrary, that there is a set $X' \subset X$ with $|X'|\geq 7$ which has all one color on edges to a large part, say blue. By Fact~\ref{Fact:Structures}, $X'$ contains neither a blue copy of $2K_2$ nor a vertex of blue degree at least $2$, and so $X'$ contains at most one blue edge. Since $X'$ consists of only small parts, the red subgraph of $X'$ must therefore be a complete graph minus a matching. Since $|X'| \geq 7$, there is a red copy of $S_6^2$ in $X'$, a contradiction.
\end{proof}

From Claim~\ref{Clm:5}, we have $|X|\leq 6$. Since $|H_{2}|\leq 5$ and $|H_{3}|\leq 5$, it follows that $|H_{1}|\geq 10$.
By Fact~\ref{Fact:Structures}, $H_1$ contains neither red copy of $2K_2$ nor a vertex of red degree at least $2$, and hence $H_1$ contains at most one red edge. If $X \neq \emptyset$, then $H_{1}$ contains no blue copy of $2K_{2}$ so by deleting at most $2$ vertices, we obtain a subgraph of $H_{1}$ of order at least $9$ with no blue edges, clearly yielding a green copy of $S_{6}^{2}$. On the other hand, if $X = \emptyset$, then $|H_{1}| \geq 16$. By deleting at most one vertex from $H_{1}$, we obtain a subgraph of $H_{1}$ of order at least $15$ with no red edges, only green and blue. By Theorem~\ref{Thm:St2Classic}, there is a monochromatic copy of $S_{6}^{2}$ in this subgraph, for a contradiction.

\begin{claim}\labelz{Claim:All7}
With $k = 3$, all parts have order at most $7$.
\end{claim}

\begin{proof}
For a contradiction, suppose $|H_{1}| \geq 8$. By minimality of $m$, there is at least one vertex in $G \setminus H_{1}$ with red edges to $H_{1}$ and at least one vertex in $G \setminus H_{1}$ with blue edges to $H_{1}$. By Fact~\ref{Fact:Structures}, $H_{1}$ contains no red or blue copy of $2K_{2}$. This means that the removal of at most $2$ vertices from $H_{1}$ can produce a subgraph with no red edges and the removal of at most an addition $2$ vertices can produce a subgraph with no blue edges. Therefore, removing at most $4$ vertices from $H_{1}$ produces a subgraph with all green edges so if $|H_{1}| \geq 10$, this subgraph contains a green copy of $S_{6}^{2}$, so $|H_{1}| \leq 9$.

With $n = 26$, there are at least $3$ vertices in $G \setminus H_{1}$ with all one color, say blue, on edges to $H_{1}$. By Fact~\ref{Fact:Structures}, there is at most one blue edge in $H_{1}$. As above, removing at most $2$ vertices from $H_{1}$ produces a subgraph with no red edges. With $|H_{1}| \geq 8$, this subgraph has order at least $6$ and all green edges except possibly one blue edge. This subgraph contains a green copy of $S_{6}^{2}$, for a contradiction.
\end{proof}

Let $A$ and $B$ be the sets of vertices in $X$ with red and blue (respectively) edges to $H_{1}$. By Claim~\ref{Clm:5}, we have $|A|, |B| \leq 6$. If $q = 1$, then by Claim~\ref{Claim:All7}, we have
\beqs
|G| & = & |H_{1}| + \sum_{i = 2}^{m} |H_{i}|\\
~ & \leq & 7 + 2 \cdot 6\\
~ & = & 19 < 26 = n,
\eeqs
a contradiction. We may therefore assume that $q = 2$. By Claims~\ref{Clm:5} and~\ref{Claim:All7}, we have
$$
|G| = |H_{1}| + |H_{2}| + |A| + |B| \leq 7 + 7 + 6 + 6 = 26 = n,
$$
which means that $|H_{1}| = |H_{2}| = 7$ and $|A| = |B| = 6$. By Claim~\ref{Clm:5}, since $|A|, |B| \geq 3$, $H_{1}$ contains at most one red edge and at most one blue edge, meaning that $|H_{1}| \leq 5$. Then
$$
|G| = |H_{1}| + |H_{2}| + |A| + |B| \leq 5 + 5 + 6 + 6 = 22 < n,
$$
a contradiction, completing the proof of Lemma~\ref{Lem:Upper-S62} when $k = 3$.




\begin{case}
$k\geq 4$.
\end{case}

As in the previous case, it is easy to see that $1 \leq q \leq 5$.

First suppose $2\leq m\leq 3$, so by minimality of $m$, we may assume $m=2$, say with corresponding parts $H_1$ and $H_2$ with all red edges between the two parts. If $|H_1|\geq 5$ and $|H_2|\geq 5$, then by Fact~\ref{Fact:Structures}, each of $H_1$ and $H_2$ contains neither a red copy of $2K_2$ nor a vertex of red degree at least $2$. Furthermore, Fact~\ref{Fact:Structures} implies that $H_1\cup H_2$ contains at most one red edge (aside from the edges in between the two parts). By deleting one vertex from $H_1\cup H_2$, there is no red edge within either part, and hence
$$
|G| = |H_1| + |H_2| \leq 2\cdot [n(k - 1) - 1] + 1 < n(k),
$$
a contradiction. Thus, suppose $|H_{2} \leq 4$ (and $|H_1|\geq 5$). Then by Fact~\ref{Fact:Structures}, $H_1$ contains no red copy of $2K_2$. By deleting two vertices from $H_{1}$, we obtain a subgraph of $H_{1}$ with no red edges, and so
$$
|G| = |H_1| + |H_2| \leq [n(k - 1) - 1] + 2 + 4 < n(k),
$$
a contradiction, meaning that we may assume $4 \leq m \leq 10$. We break the remainder of the proof into subcases based on the value of $q$.

If $q=5$, then by Fact~\ref{Fact:Structures}, $m = 5$ so $G = H_1\cup H_2\cup H_3\cup H_4\cup H_5$ and the reduced graph is the unique $2$-colored copy of $K_{5}$ with no monochromatic triangle. By Fact~\ref{Fact:Structures}, aside from edges between the parts, $H_1\cup H_2\cup H_3\cup H_4\cup H_5$ contains at most one red edge and one blue edge. Thus, by deleting at most $2$ vertices, there is neither a red nor a blue edge within any part of the Gallai partition. This implies that
$$
|G|=\sum_{i=1}^5|H_i| \leq 5 \cdot [n(k - 2) - 1] + 2 < n(k),
$$
a contradiction.

Next suppose $q = 4$. In order to avoid a monochromatic triangle among the large parts, there are only two possible colorings for the reduced graph restricted to the large parts, one with complementary monochromatic copies of $P_{4}$ and the other with a $C_{4}$ in one color and $2K_{2}$ in the other. In either case, each large part has red and blue edges to another large part so by Fact~\ref{Fact:Structures}, there are at most $2$ total red edges and at most $2$ total blue edges within large parts. Then by removing at most $4$ total vertices, we may obtain a subgraph in which there are no red or blue edges within the parts. If $m = 4$, then
$$
|G|=\sum_{i=1}^4|H_i| \leq 4[n(k - 2) - 1] + 4 < n(k),
$$
a contradiction, meaning that we may assume $m \geq 5$. If the reduced graph restricted to the large parts contains a $C_{4}$ in one color, then the edges from $X$ to the large parts cannot be colored without creating a monochromatic copy of $S_{6}^{2}$, so we may assume the reduced graph restricted to the large parts contains two complementary monochromatic copies of $P_{4}$. In order to avoid a monochromatic triangle in the reduced graph using at least two large parts, the edges from all of $X$ to the ends of each $P_{4}$ must have the same color as the $P_{4}$ itself. By the minimality of $m$, $X$ must be a single part of the Gallai partition of order at most $2$. By Fact~\ref{Fact:Structures} (as in the proof of the case $q = 5$ above), each part $H_{i}$ for $1 \leq i \leq 4$ contains at most one red edge and at most one blue edge. Then by removing at most $2$ total vertices, we may obtain a subgraph in which there are no red or blue edges within the parts. Thus,
$$
|G|=\sum_{i=1}^4|H_i| + |X| \leq 4[n(k - 2) - 1] + 2 + 2 < n(k),
$$
a contradiction.


Suppose $q=3$. Suppose that the edges from $H_1$ to $H_2\cup H_3$ are red and the edges from $H_2$ to $H_3$ are blue. Then to avoid a monochromatic triangle using at least two large parts, every vertex in $X$ has blue edges to $H_{1}$. By Claim~\ref{Clm:5}, we have $|X| \leq 6$. By Fact~\ref{Fact:Structures}, $H_1\cup H_2\cup H_3$ contains at most one red edge within the parts, $H_2\cup H_3$ contains at most one blue edge within the parts, and $H_1$ does not contain blue copy of $2K_2$. By deleting at most $4$ vertices, we may obtain a subgraph in which there is neither red nor blue edges within the parts of $H_1\cup H_2\cup H_3$, and hence
$$
|G| = \sum_{i=1}^3 |H_i| + |X| \leq 3[n(k - 2) - 1] + 4 + 6 < n(k),
$$
a contradiction.

Suppose $q=2$, say with the edges from $H_1$ to $H_2$ being red. By Fact~\ref{Fact:Structures}, $H_1\cup H_2$ contains at most one red edge within the parts, and each large part $H_i$ does not contain a blue copy of $2K_2$. By deleting at most $5$ total vertices from the large parts, we may obtain a subgraph in which neither red nor blue edges appear within the parts of $H_1\cup H_2$. Let $A$ be the set of small parts with red edges to $H_1$, and let $B$ be the set of small parts with blue edges to $H_1$. By Claim~\ref{Clm:5}, we have $|A|\leq 6$ and $|B|\leq 6$, and hence
$$
|G| = |H_1| + |H_2| + |A| + |B| \leq 2[n(k - 2) - 1] + 5 + 12 < n(k),
$$
a contradiction.

Finally suppose $q=1$ and let $A$ be the set of small parts with red edges to $H_1$ and $B$ be the set of small parts with blue edges to $H_1$. From Claim~\ref{Clm:5}, we have $|A|\leq 8$ and $|B|\leq 8$. By Fact~\ref{Fact:Structures}, $H_1$ contains neither a red copy of $2K_2$ nor a blue copy of $2K_2$. By deleting at most $4$ vertices from $H_{1}$, we may obtain a subgraph of $H_{1}$ in which there is neither a red nor a blue edge. Hence
$$
|G| = |H_1| + |A| + |B| \leq [n(k - 2) - 1] + 4 + 12 < n(k),
$$
a contradiction.
\end{proof}

\subsection{The case $r=2$ and $t=8$}

We first give a lower bound construction for Item $(2)$ of Theorem~\ref{Thm:St2}.

\begin{lemma}\labelz{Lem:Lower-S82}
For $k \geq 3$,
$$
gr_k(K_3;S_8^{2})\geq \begin{cases}
14 \cdot 5^{\frac{k-2}{2}} + \frac{1}{2} \cdot 5^{\frac{k - 4}{2}} + \frac{1}{2},  & \text{ if $k$ is even;}\\
7\cdot 5^{\frac{k-1}{2}}+\frac{1}{4}\cdot 5^{\frac{k-3}{2}}+\frac{3}{4}, & \text{ if $k$ is odd.}
\end{cases}
$$
\end{lemma}

\begin{proof}
This lower bound result is proven by constructing a colored complete graph on one fewer vertex that does not have the desired colored subgraphs.

Let $G_1$ be a copy of $K_7$ entirely colored by color $1$ and let $G_2$ be a copy of $K_{14}$ consisting of two copies of $G_1$ joined by all edges of color $2$. Now suppose we have constructed a Gallai colored complete graph $G_{k-2}$ with $k \geq 3$ using $k-2$ colors which contains no monochromatic copy of $S_8^2$. We construct the graph $G_k$ by making five copies of $G_{k-2}$ and inserting edges of colors $k-1$ and $k$ between the copies to form a blow-up of the unique $2$-coloring of $K_5$ with no monochromatic triangle. This coloring contains no rainbow triangle and no monochromatic copy of $S_8^2$. When $k = 3$, we set $G'_3 = G_3$ and $G'_{4} = G_{4}$.

For each odd $k \geq 5$, we extend this construction as follows. For $j \geq 5$, let $A_1^j$ be a colored copy of $K_8$ containing a perfect matching with colors $2, 3, j-1$ and $j$ and all other edges with color $1$. Now choose one copy of $G_{1}$ from the construction of $G_{k}$ and replace it with a copy of $A_1^k$ to produce a slightly larger graph $G'_{k}$.

For each even $k \geq 6$, we extend this construction as follows. For $j \geq 6$, let $A_{1}^{j}$ be a colored copy of $K_{8}$ containing a perfect matching with colors $2, 3, 4, j$ and all other edges having color $1$. Similarly let $A_{2}^{j}$ be a colored copy of $K_{8}$ containing a perfect matching with colors $2, 3, 4, j - 1$ and all other edges having color $1$. If $k = 6$, we choose two of the copies of $G_{4}$ and choose one copy of $G_{1}$ from each. One of these copies of $G_{1}$ is replaced by $A_{1}^{6}$ and the other is replaced by $A_{2}^{6}$ to create the graph $G'_{6}$ on $352$ vertices. For $k \geq 8$, there are always at least $15$ copies of $G_{4}$ (which had not been modified) so choose two of these and again replace two copies of $G_{1}$ by $A_{1}^{k}$ and $A_{2}^{k}$ respectively.

This modification process is applied at each step of the construction process so in particular, for odd values of $k$, we have
\bi
\item $|G'_{3}| = 35$
\item $|G'_{5}| = 5 \cdot 35 + 1 = 176$,
\item $|G'_{7}| = 5 \cdot 176 + 1 = 881$.
\ei
Similarly, for even values of $k$, we have
\bi
\item $|G'_{4}| = 5 \cdot 14 = 70$
\item $|G'_{6}| = 5 \cdot 70 + 1 = 352$
\item $|G'_{8}| = 5 \cdot 352 + 2 = 1762$
\ei
If $k$ is odd and $k \geq 5$, this yields
\beqs
|G'_k| & = & 5 \cdot |G'_{k - 2}| + 1\\
~ & = & 5\left(7\cdot 5^{\frac{k-3}{2}}+\frac{1}{4}\cdot 5^{\frac{k-5}{2}}+\frac{3}{4} - 1\right) + 1\\
~ & = & 7\cdot 5^{\frac{k-1}{2}}+\frac{1}{4}\cdot 5^{\frac{k-3}{2}} - \frac{1}{4}.
\eeqs
Similarly, if $k$ is even and $k \geq 6$, this yields a new graph $G'_{k}$ with
\beqs
|G'_k| & = & 5 \cdot |G'_{k - 2}| + 2\\
~ & = & 5 \left( 14 \cdot 5^{\frac{k-4}{2}} + \frac{1}{2} \cdot 5^{\frac{k - 6}{2}} + \frac{1}{2} - 1 \right) + 2\\
~ & = & 14 \cdot 5^{\frac{k-2}{2}} + \frac{1}{2} \cdot 5^{\frac{k - 4}{2}} - \frac{1}{2},
\eeqs
as claimed.
\end{proof}

We now give an upper bound for $(1)$ of Theorem \ref{Thm:St2}.

\begin{lemma}\labelz{Lem:S82}
For $k \geq 1$,
$$
gr_k(K_3;S_8^{2})\leq \begin{cases}
8, & \text{ if $k = 1$;}\\
15, & \text{ if $k = 2$;}\\
14\cdot 5^{\frac{k-2}{2}}+\frac{1}{4}\cdot 5^{\frac{k-4}{2}}+\frac{3}{4}, & \text{ if $k$ is even;}\\
7\cdot 5^{\frac{k-1}{2}}+\frac{1}{4}\cdot 5^{\frac{k-3}{2}}+\frac{3}{4}, &\text{ if $k$ is odd.}
\end{cases}
$$
\end{lemma}

\begin{proof}
We prove the upper bound by induction on $k$. The case $k=1$ is trivial and the case $k=2$ is precisely $R(S_8^{2},S_8^{2})=15$. We therefore suppose $k\geq 3$ and let $G$ be a coloring of $K_n$ with no rainbow triangle and no monochromatic copy of $S_8^{2}$ where
$$
n= n_{k} = \begin{cases}
14\cdot 5^{\frac{k-2}{2}}+\frac{1}{4}\cdot 5^{\frac{k-4}{2}}+\frac{3}{4}, &\text{ if $k$ is even;}\\
7\cdot 5^{\frac{k-1}{2}}+\frac{1}{4}\cdot 5^{\frac{k-3}{2}}+\frac{3}{4}, &\text{ if $k$ is odd.}
\end{cases}
$$

Since $G$ is a Gallai coloring, it follows from Theorem~\ref{Thm:G-Part} that there is a Gallai partition of $V(G)$. Suppose that the two colors appearing in the Gallai partition are red and blue. Let $m$ be the number of parts in this partition and choose such a partition where $m$ is minimized. Let $H_{1}, H_{2}, \dots, H_{m}$ be the parts of this partition, say with $|H_{1}| \geq |H_{2}| \geq \dots \geq |H_{m}|$.
Let $q$ be the number of ``large'' parts of order at least $4$. This means that $|H_{q}| \geq 4$ and $|H_{q + 1}| \leq 3$. Let $X$ be the set of vertices in the ``small'' parts of order at most $3$. Note that $m \leq 14$ by Theorem~\ref{Thm:St2Classic}. First some easy observations much like Fact~\ref{Fact:Structures}.

\begin{fact}\labelz{Fact:S82Struc}
Let $A$ and $B$ be two disjoint non-empty sets of vertices in a colored complete graph with no monochromatic copy of $S_{8}^{2}$, with all edges from $A$ to $B$ colored red.
\bi
\item If $|B| \geq 5$, then $A$ contains no vertex with two incident red edges (within $A$).
\item If $|A| \geq 7$, then $A$ contains no red copy of $2K_{2}$.
\item If $|A| \geq 6$, then at most one of $A$ or $B$ contains a red edge.
\ei
Furthermore, given three disjoint non-empty sets of vertices with all red edges between pairs of these sets, then either the two larger sets have total order at most $6$ or two of the sets have order $1$. In particular, there can be no monochromatic triangle in the reduced graph that corresponds to at least two large parts (of order at least $4$).
\end{fact}

We break the proof into two cases, when $k = 3$ and when $k \geq 4$. 

\setcounter{case}{0}
\begin{case}
$k = 3$ so $n = 36$.
\end{case}

With red and blue being the colors appearing on edges between parts of the Gallai partition, let green be the third color in this $3$-coloring.

First suppose $2\leq m\leq 3$ so by the minimality of $m$, we may assume $m = 2$, say with all red edges between $H_{1}$ and $H_{2}$. Clearly, $|H_1|\geq \frac{n}{2} = 18$ so by Fact~\ref{Fact:S82Struc}, $H_{1}$ contains no red copy of $2K_{2}$. By removing at most $2$ vertices from $H_{1}$, we obtain a subgraph with no red edges. Since $|H_{1}| - 2 \geq 16$ and $R(S_8^{2},S_8^{2})=15$, it follows that there is a monochromatic copy of $S_8^{2}$ within $H_{1}$, a contradiction. We may therefore assume that $4 \leq m \leq 14$. By the minimality of $m$, this means every part has incident edges (to other parts) in both red and blue.

The remainder of this case is broken into subcases based on the value of $q$.

First suppose $q = 5$. Then to avoid a monochromatic triangle in the reduced graph containing at least two large parts, we must also have $m = 5$ and the reduced graph must be the unique $2$-coloring of $K_{5}$ with no monochromatic triangle. Then $|H_{1}| \geq \ceil{\frac{36}{5}} = 8$. By Fact~\ref{Fact:S82Struc}, there is a total of at most one red edge and at most one blue edge inside the parts $H_1\cup H_2\cup H_3\cup H_4\cup H_5$. Hence, $H_{1}$ contains at most one red edge and at most one blue edge, with all remaining edges being green. This produces a green copy of $S_{8}^{2}$, a contradiction.

Next suppose $q = 0$. By the Pigeonhole Principle, we have $|H_{1}| \geq 3$ so let $A$ be the set of vertices in parts with red edges to $H_{1}$ and let $B$ be the set of vertices in parts with blue edges to $H_{1}$. Since $|A \cup B| = 33$, it follows that one of $|A|\geq 17$ or $|B|\geq 17$ must hold, say $|A|\geq 17$ without loss of generality. By Fact~\ref{Fact:S82Struc}, there can be no red copy of $2K_{2}$ within $A$. By deleting at most $2$ vertices from $A$, we may obtain a subgraph with no red edges. Since $|A|-2\geq 15 = R(S_8^{2},S_8^{2})$, it follows that this subgraph of $A$ contains either a blue or a green copy of $S_8^{2}$, a contradiction.

Next suppose $q = 1$. By Fact~\ref{Fact:S82Struc}, if $|H_{1}| \geq 7$, there can be no red or blue copy of $2K_{2}$ within $H_{1}$. By removing at most $2$ vertices from $H_{1}$, we obtain a subgraph in which there is at most one red edge and at most one blue edge. If $|H_{1}| \geq 10$, such a subgraph would contain a green copy of $S_{8}^{2}$, so we must have $|H_{1}| \leq 9$. Let $A$ be the set of vertices in parts with red edges to $H_1$ and let $B$ be the set of vertices in parts with blue edges to $H_1$. Since $n=36$, it follows that $|A|+|B|\geq 27$, so either $|A|\geq 14$ or $|B|\geq 14$, say $|A|\geq 14$. By Fact~\ref{Fact:S82Struc}, $A$ contains no red copy of $2K_2$. This means that by deleting at most $2$ vertices from $A$, we have a subgraph with no red edges. Since $|A| - 2 \geq 12$ and all parts of $A$ have order at most $3$, there are at least $4$ parts within this subgraph with all blue edges in between each pair. This produces a blue copy of $S_8^{2}$, for a contradiction.

Next suppose $q=4$. If $m \geq 5$, then in order to avoid a monochromatic triangle in the reduced graph including two large parts, all small parts must have specifically colored edges to the large parts such that the small parts together form one part of a blow-up of the unique $2$-coloring of $K_{5}$ with no monochromatic triangle. By the minimality of $m$, there is only one part $H_{5}$ other than $H_1\cup H_2\cup H_3\cup H_4$, with $|H_5|\leq 3$. Since $n=36$, regardless of the value of $m$, it follows that $\sum_{i=1}^4|H_i|\geq 33$. By the Pigeonhole Principle, we have $|H_{1}| \geq \ceil{\frac{33}{4}} = 9$ and by Fact~\ref{Fact:S82Struc}, there can be no red or blue copy of $2K_{2}$ within $H_{1}$. By removing at most $2$ vertices from $H_{1}$, we obtain a subgraph in which there is at most $1$ red edge and at most one blue edge. If $|H_{1}| \geq 10$, such a subgraph would contain a green copy of $S_{8}^{2}$, so we must have $|H_{1}| = 9$. Since $\sum_{i=1}^4|H_i|\geq 33$, we have $|H_{2}| \geq \ceil{\frac{24}{3}} = 8$ and suppose the edges from $H_{1}$ to $H_{2}$ are red. By Fact~\ref{Fact:S82Struc}, $H_{1}$ contains at most one red edge. By removing at most one vertex from $H_{1}$, we obtain a subgraph in which there is at most one blue edge. Since $|H_{1}| - 1 \geq 8$, this subgraph contains a green copy of $S_{8}^{2}$, for a contradiction.

Next suppose $q = 3$. Disregarding the relative orders of the large parts, we may assume without loss of generality and to avoid a monochromatic triangle, that all edges from $H_{1}$ to $H_{2} \cup H_{3}$ are red and the edges from $H_{2}$ to $H_{3}$ are blue. Then in order to avoid a monochromatic triangle in the reduced graph corresponding to two of the large parts, all edges from $X$ to $H_{1}$ must be blue. Next we claim that the large parts are not too large.

\begin{claim}\labelz{Clm:8}
$|H_i|\leq 7$ for all $i$.
\end{claim}

\begin{proof}
First suppose $|H_2|\geq 8$. By Fact~\ref{Fact:S82Struc}, $H_2$ contains neither a red nor a blue copy of $2K_2$. Furthermore, since $H_{1}$ and $H_{3}$ are both large, $H_2$ contains neither a vertex of red degree at least $3$ nor a vertex of blue degree at least $3$ to avoid creating a monochromatic copy of $S_{8}^{2}$. Therefore, the graph induced by the red edges in $H_2$ is a subgraph of a red triangle, and the graph induced by the blue edges in $H_2$ is a subgraph of a blue triangle. The remaining edges of $H_{2}$ all being green, there is a green copy of $S_8^{2}$ in $H_{2}$, a contradiction. This means that $|H_2| \leq 7$ and symmetrically $|H_{3}| \leq 7$.

Next assume $|H_1|\geq 8$. By Fact~\ref{Fact:S82Struc}, $H_1$ contains at most one red edge, and no blue copy of $2K_2$. Furthermore, the blue edges within $H_1$ form either a blue triangle or a blue star. If the blue edges of $H_{1}$ are a subgraph of a triangle, then with all remaining edges being green aside from possibly one red edge, there is a green copy of $S_8^{2}$ in $H_{1}$, a contradiction. Thus suppose that $H_1$ contains a blue star with at least $3$ edges. By Fact~\ref{Fact:S82Struc}, we must have $|X|\leq 4$. With $|H_{2}| + |H_{3}| + |X| \leq 7 + 7 + 4 = 18$ and $n = 36$, we have $|H_{1}| \geq 18$. By removing at most two vertices from $H_{1}$, we obtain a subgraph in which there are no red or blue edges. Since this subgraph has order at least $16$ and all green edges, it must contain a green copy of $S_{8}^{2}$ for a contradiction.
\end{proof}

From Claim~\ref{Clm:8}, we have $|H_i|\leq 7$ for $i=1,2,3$. Since $n=36$, it follows that $|X|\geq 15$. By Fact~\ref{Fact:S82Struc}, there is no blue copy of $2K_{2}$ in $X$ so by deleting at most $2$ vertices from $X$, we can obtain a subgraph of $X$ of order at least $13$ in which there are no blue edges. Since $X$ consists of only small parts of the Gallai partition, at least $5$ of them, and these parts must have all red edges in between them, it follows that there is a red copy of $S_{8}^{2}$ as a subgraph of $X$, completing the proof when $q = 3$.

Finally suppose $q=2$. Suppose that the edges from $H_1$ to $H_2$ are red. Let $A$ be the set of vertices in small parts with red edges to $H_1$, and let $B$ be the set of vertices in small parts with blue edges to $H_1$. Then we have the following claim.

\begin{claim}\labelz{Clm:9}
$|H_i|\leq 9$ for all $i$.
\end{claim}

\begin{proof}
Assume, to the contrary, that $|H_1|\geq 10$. By Fact~\ref{Fact:S82Struc}, $H_1$ contains neither a red nor a blue copy of $2K_2$. By deleting at most $2$ vertices from $H_{1}$, we can obtain a subgraph in which there is at most one red edge and at most one blue edge. Since all other edges of this subgraph are green and the subgraph has order at least $|H_{1}| - 2 \geq 8$, this contains a green copy of $S_{8}^{2}$, a contradiction.
\end{proof}

From Claim~\ref{Clm:9}, we have $|H_1|, |H_{2}| \leq 9$. Let $v \in H_{2}$ and let $A' = A \cup \{v\}$. By Fact~\ref{Fact:S82Struc}, $A'$ contains no red copy of $2K_{2}$ so removing at most $2$ vertices from $A'$ yields a subgraph in which there are no red edges. If $|A| \geq 9$, then such a subgraph of $A'$ has order at least $8$ and consists of only small parts of the Gallai partition, in between which are all blue edges, making a blue copy of $S_{8}^{2}$. This means that $|A| \leq 8$. By the same argument, we find that $|B| \leq 9$. This is a contradiction to the fact that $n = 36$ since
$$
n = |H_{1}| + |H_{2}| + |A| + |B| \leq 9 + 9 + 8 + 9 = 35.
$$

\begin{case}
$k \geq 4$.
\end{case}

First suppose $2\leq m\leq 3$, so we may assume $m = 2$. If $|H_1|\geq |H_2|\geq 7$, then the parts $H_1, H_2$ together contain at most one red edge, and hence $|G|=|H_1|+|H_2|\leq 2(n_{k - 1}-1)+1 < n_{k}$, a contradiction. If $|H_1|\geq 7$ and $|H_2|\leq 6$, then $H_1$ contains at most one red edge, and hence $|G|=|H_1|+|H_2|\leq n_{k - 1}+6 < n_{k}$, a contradiction. If $|H_1|\leq 6$ and $|H_2|\leq 6$, then $|G|=|H_1|+|H_2|\leq 12 < n_{k}$, a contradiction for all $k \geq 4$.

Since $R(S_8^{2},S_8^{2})=15$, it follows that $4\leq m\leq 14$. Since $k \geq 4$, we have $n \geq 71$ so by the pigeonhole principle, there is a part of order at least $\ceil{ \frac{71}{14} } = 6$. Let $a$ be the positive integer such that $|H_i|\geq 4$ for $1\leq i\leq a$, and $|H_i|\leq 3$ for $a+1\leq i\leq m$. Note that $1 \leq a \leq 5$ since there is at least one part of order at least $6$ and if there were at least $6$ such parts, there would be a triangle in the reduced graph among those parts of order at least $4$, making a monochromatic copy of $S_{8}^{2}$. Let $A = \bigcup_{i > a} H_{i}$.

If $a=1$, then since $m\leq 14$ and $n \geq 71$, we have $|H_1|\geq 71 - 13\cdot 3 = 32$. By Fact~\ref{Fact:S82Struc}, we know that $H_1$ contains neither a red copy of $2K_2$ nor a blue copy of $2K_2$. By deleting a total of at most $4$ vertices, we can obtain a subgraph of $H_{1}$ in which there are no red or blue edges. We therefore have
$$
|G|=\sum_{i=1}^m|H_i|\leq (n_{k - 2}-1)+4+3(m-1) \leq n_{k - 2}+42<n_{k},
$$
a contradiction.

Suppose $a=2$. If $|H_1|\geq 7$ and $|H_2|\geq 7$, then by Fact~\ref{Fact:S82Struc}, $H_1\cup H_2$ contains at most one red edge, $H_1$ contains no blue $2K_2$, and $H_2$ contains no blue $2K_2$. By deleting at most $5$ vertices, we can obtain a subgraph of $H_{1} \cup H_{2}$ in which there is no red and blue edge within either $H_1$ or $H_2$. This means that 
$$ 
|G| = \sum_{i=1}^m|H_i| \leq 2(n_{k - 2} -1) + 5 + 3(m-2) = 2(n_{k - 2} -1) + 41 < n_{k}
$$
for all $k \geq 4$, a contradiction.

Next suppose $a=3$. In order to avoid a monochromatic copy of $S_{8}^{2}$, the reduced graph restricted to the vertices corresponding to $\{H_{1}, H_{2}, H_{3}\}$ must not form a monochromatic triangle. Therefore, suppose that the edges from $H_1$ to $H_2\cup H_3$ are red, the edges from $H_2$ to $H_3$ are blue. If there is a part $H_{i}$ with red edges to $H_{1}$, then it must have blue edges to both $H_{2}$ and $H_{3}$, creating a blue copy of $S_{8}^{2}$. Thus to avoid a monochromatic copy of $S_8^2$, for each $i \geq 4$, the edges from $H_1$ to $H_i$ must be blue. By Fact~\ref{Fact:S82Struc}, $H_{2} \cup H_{3}$ contains no red copy of $2K_{2}$ so the removal of at most $2$ vertices leaves behind a subgraph of $H_{2} \cup H_{3}$ with no red edges. If $|H_{i}| \geq 7$ for $i \in \{2, 3\}$, then by Fact~\ref{Fact:S82Struc}, there is a total of at most one blue edge within $H_{2}$ and $H_{3}$. Thus, by removing at most one additional vertex, we obtain a subgraph of $H_{2} \cup H_{3}$ in which there are no blue edges within the parts. We therefore have
$$
|H_2|+|H_3| \leq 2(n_{k - 2} - 1) + 3.
$$

First if $|H_1|\geq 7$. By Fact~\ref{Fact:S82Struc}, $H_{1}$ contains no red or blue copy of $2K_{2}$ and since $|H_{2} \cup H_{3}| \geq 8$, there is also no vertex with two incident red edges within $H_{1}$. This means that by removing at most $3$ vertices from $H_{1}$, we can destroy all red and blue edges within $H_{1}$ so $|H_{1}| \leq (n_{k - 2} - 1) + 3$. Also by Fact~\ref{Fact:S82Struc}, $A$ contains at most one blue edge. Then $|A|\leq 15$ since otherwise $A$ consists of at least $6$ parts, thereby containing a red copy of $S_{8}^{2}$. This means we have 
$$
|G|=\sum_{i=1}^3|H_i|+|A|\leq 3(n_{k - 2} - 1) + 6 + 15 = 3n_{k - 2} + 18 < n_{k},
$$
a contradiction.

Next if $|H_1|\leq 6$, then we claim that $|A|\leq 17$. Indeed, if $|A|\geq 18$, then since $A$ contains no blue $2K_2$, the blue edges in $A$ form a triangle or a star. By deleting $2$ vertices, there is no blue edge within $A$. Since $|A|-2\geq 16$, it follows that there are at least six parts with all red edges in between, providing a red copy of $S_{8}^{2}$. We therefore have 
$$
|G|=\sum_{i=1}^{3} |H_i| + |A| \leq 6 + 2(n_{k - 2} - 1) + 3 + 17 = 2n_{k - 2} + 24 < n_{k},
$$
a contradiction.

Next suppose $a=4$. If $m\geq 5$, then for all $i \leq 4$, all edges from $H_{i}$ to $A$ have a single color and the reduced graph contains a red cycle $H_1H_2H_3H_4AH_1$ and a blue cycle $H_1H_3AH_2H_4H_1$. By Fact~\ref{Fact:S82Struc}, there is no red copy of $2K_{2}$ within $H_{1} \cup H_{4}$ and if either of $H_{1}$ or $H_{4}$ has order at least $7$, then it contains no blue copy of $2K_{2}$. Also if either $H_{1}$ or $H_{2}$ has order at least $6$, then at most one of the two sets can contain any blue edges. Putting these observations together, we have $|H_1|+|H_4|\leq 2(n_{k - 2}-1)+4$ and similarly $|H_2|+|H_3|\leq 2(n_{k - 2} -1)+4$. By the choice of the Gallai partition with $m$ minimum, we have $A = H_{5}$ so $|A|\leq 3$. Then 
$$
|G|=\sum_{i=1}^5 |H_i| \leq 4(n_{k - 2}-1) + 8 + 3 < n_{k},
$$
a contradiction.  If $m=a=4$, then again we have $|H_1|+|H_4|\leq 2(n_{k - 2}-1)+4$ and $|H_2|+|H_3|\leq 2(n_{k - 2} -1)+4$ so 
$$
|G|=\sum_{i=1}^4 |H_i| \leq 4(n_{k - 2}-1) + 8 < n_{k},
$$
a contradiction.

Finally suppose $a = m = 5$. By Fact~\ref{Fact:S82Struc}, we have that $H_1\cup H_2\cup H_3\cup H_4\cup H_5$ contains a total of at most one red edge and at most one blue edge within the parts. If $k\geq 5$, then $|G|=\sum_{i=1}^5|H_i|\leq 5(n_{k - 2}-1)+2 <n_{k}$, a contradiction. We may therefore assume that $k = 4$ so $n = 71$.

Since $\ceil{ \frac{71}{5} } = 15$, there exists a part, say $H_{1}$, such that $|H_{1}|\geq 15$. Since $R(S_{8}^{2}) = 15$, there is at least one red or blue edge within $H_{1}$. Without loss of generality, suppose $w_{1}w_{2}$ is a red edge within $H_{1}$ and if there is a blue edge within $H_{1}$, let $w_{3}w_{4}$ be such an edge. Note that these two edges do not share a vertex to avoid a rainbow triangle. All other edges of $H_{1}$ are colored by colors $3$ and $4$, say green and purple. Let $F = H_{1} \setminus \{w_{1}, w_{2}, w_{3}, w_{4}\}$ so every vertex in $F$ has only green or purple edges to the rest of $H_{1}$.

At this point, it is worth noting that all other parts of the Gallai partition may have order at most $14$ and so may be $2$-colored using green and purple to avoid a monochromatic copy of $S_{8}^{2}$. With $H_{1}$ containing at most one red and at most one blue edge, the remainder of the proof consists of finding either a rainbow triangle or a green or purple copy of $S_{8}^{2}$ within $H_{1}$. We may also therefore assume that $|H_{1}| = 15$.

Since $|H_{1}| = 15$, there are at least $7$ vertices in $H_{1} \setminus \{w_{1}, w_{2}\}$ with all one color on their edges to $\{w_{1}, w_{2}\}$. Without loss of generality, say this color is purple and let $P$ be the set of vertices in $H_{1}$ with all purple edges to $\{w_{1}, w_{2}\}$. By Fact~\ref{Fact:S82Struc}, $P$ contains no purple copy of $2K_{2}$. This means that the purple edges within $P$ form either a star or a triangle so almost all edges within $P$ are green. Thus, if $|P| \geq 9$, then $P$ contains a green copy of $S_{8}^{2}$ so $7 \leq |P| \leq 8$. Let $Q = H_{1} \setminus (P \cup \{w_{1}, w_{2}\})$.

First suppose $|P| = 8$ so $|Q| = 5$. In order to avoid a green copy of $S_{8}^{2}$ within $P$, there must exist a vertex $v \in P$ that is the center of a spanning purple star of $P$. Then with $w_{1}$ and $w_{2}$ each forming a purple triangle, $v$ is the center of a purple copy of $S_{8}^{2}$, for a contradiction.

Thus, we may assume that $|P| = 7$ so $|Q| = 6$. First suppose $P$ contains a purple star with at least $3$ edges, say centered at $v$, so no purple triangle. Then every vertex of $P$ except $v$ has at most one incident purple edge within $P$. To avoid a green copy of $S_{8}^{2}$, each vertex of $P \setminus \{v\}$ has at most one green edge to $Q$. Thus, every vertex of $P \setminus \{v\}$ has at least $5$ purple edges to $Q$. 
Let $u \in P \setminus \{v\}$ with a purple edge to $v$ and let $u'$ be the vertex of $Q$ (if one exists) with a green edge to $u$. Then $Q \setminus \{u'\}$ contains no purple edge to avoid a purple copy of $S_{8}^{2}$. To avoid a purple copy of $S_{8}^{2}$, $v$ can have at most $4$ purple edges to $Q$, leaving at least two edges to $Q$ which must be green. At least one of these must go to a vertex $v' \in Q \setminus \{u'\}$, forming a green copy of $S_{8}^{2}$ centered at $v'$.

Therefore, we may assume $P$ contains no vertex with purple degree at least $3$, or rather, the purple edges within $P$ form a subgraph of a triangle, say $T$. In order to avoid a green copy of $S_{8}^{2}$, every vertex of $P \setminus T$ has all purple edges to $Q$. Since $|Q \cup \{w_{1}, w_{2}\}| = 8$, by Fact~\ref{Fact:S82Struc}, there is no purple copy of $2K_{2}$ within $Q$. To avoid a green copy of $S_{8}^{2}$ centered within $Q$, there must be a vertex $x \in Q$ with purple edges to all of $Q \setminus \{x\}$. Then using two vertices of $P$ to form purple triangles, $x$ is the center of a purple copy of $S_{8}^{2}$, to complete the proof.
\end{proof}



\newpage

\appendix

\section{Appendix for review}

\subsection{The case for $r=2$ and general $t$}

In this section, we prove Item~$(3)$ of Theorem~\ref{Thm:St2}. First we prove the following lemma, which provides the lower bound.

\begin{lemma}\labelz{Lem:Lower-St2}
$$
gr_k(K_3;S_t^{2})\geq \begin{cases}
2(t-1)\times 5^{\frac{k-2}{2}}+1, & \text{ if $k$ is even;}\\
(t-1)\times 5^{\frac{k-1}{2}}+1, &\text{ if $k$ is odd.}
\end{cases}
$$
\end{lemma}

\begin{proof}
We prove this result by inductively constructing a coloring of $K_{n}$ where
$$
n=\begin{cases}
2(t-1)\cdot5^{\frac{k-2}{2}} & \text{ if $k$ is even,}\\
(t-1)\cdot5^{\frac{k-1}{2}} & \text{ if $k$ is odd,}
\end{cases}
$$
which contains no rainbow triangle and no monochromatic copy of $S_{t}^{2}$. 

If $k$ is odd, let $G_{1}$ be a $1$-colored complete graph on $t-1$ vertices. Without sufficient vertices, this contains no monochromatic copy of $S_{t}^{2}$. Suppose this coloring uses color $1$. Suppose we have constructed a coloring $G_{2i-1}$ where $i$ is a positive integer and $2i-1<k$, using the $2i-1$ colors $1,2, \dots, 2i-1$ and having order $n_{2i-1}=(t-1)\cdot5^{i-1}$. Construct $G_{2i+1}$ by making five copies of $G_{2i-1}$ and inserting edges of color $2i$ and $2i+1$ between the copies to form a blow-up of the unique $2$-colored $K_{5}$ which contains no monochromatic triangle. This coloring clearly contains no rainbow triangle and, since there is no monochromatic triangle in either of the two new colors, there can be no monochromatic copy of $S_{t}^{2}$ in $G_{2i+1}$.

If $k$ is even, let $G_{k - 1}$ be as constructed in the odd case above. Construct $G_{k}$ by making two copies of $G_{k - 1}$ and inserting all edges of color $k$ in between the two copies. This graph certainly contains no rainbow triangle and since color $k$ is bipartite, it also contains no monochromatic copy of $S_{t}^{2}$, and has order 
$$
|G_{k}| = 2|G_{k - 1}| = 2(t - 1) \cdot 5^{\frac{k - 2}{2}},
$$
as desired.
\end{proof}

We are now in a position to complete the proof of Item~$(3)$ from Theorem~\ref{Thm:St2}, that is, to prove that
$$
gr_k(K_3,S_{t}^{2})\leq \begin{cases}
2t\cdot 5^{\frac{k-2}{2}}, &\mbox {\rm if}~k~is~even;\\[0.2cm]
t\cdot 5^{\frac{k-1}{2}}, &\mbox {\rm if}~k~is~odd.
\end{cases}
$$

\noindent {\it Proof of Item~$(3)$ of Theorem~\ref{Thm:St2}}. 
The lower bound follows from Lemma~\ref{Lem:Lower-St2}. We prove the upper bound by induction on $k$. The case $k=1$ is immediate and the case $k=2$ is precisely the result of Theorem~\ref{Thm:St2Classic}, so suppose $k\geq3$ and let $G$ be a coloring of $K_{n}$ where
$$
n=\begin{cases}
2t\cdot5^{\frac{k-2}{2}} & \text{ if $k$ is even,}\\
t\cdot5^{\frac{k-1}{2}}& \text{ if $k$ is odd.}
\end{cases}
$$
For a contradiction, we suppose $G$ is a Gallai coloring which contains no monochromatic copy of $S_{t}^{2}$.

By Theorem~\ref{Thm:G-Part}, there is a Gallai partition of $G$ and suppose red and blue are the two colors appearing in the Gallai partition. Let $m$ be the number of parts in this partition and choose such a partition where $m$ is minimized. By Theorem~\ref{Thm:G-Part}, since choosing one vertex from each part of the partition yields a $2$-colored complete graph, we see that $m\leq 2t-2$. Let $r$ be the number of ``large'' parts of the Gallai partition with order at least $t$, say with $|H_{1}|,|H_{2}|, \ldots, |H_{r}|\geq t$ and $|H_{r+1}|,|H_{r+2}|,\ldots, |H_{m}|\leq t-1$. Call all remaining parts, those with order at most $t - 1$, ``small''.

First some helpful supporting observations. 

\begin{claim}\labelz{Claim:Claim1}
Suppose $m \geq 3$.
\begin{itemize}
\item $(1)$ If $r = 2$ and $A$ is the set of vertices with all blue edges to $H_{1}$ and all red edges to $H_{2}$ (or equivalently red and blue respectively), then $|A|\leq t-1$.
\item $(2)$ If $r \geq 1$ and $B$ is a set of vertices in small parts (of order at most $t - 1$) with all one color on edges to $H_{1}$, then $|B|\leq 2(t-1)$.
\item $(3)$ If $B$ is a set of vertices in small parts with all one color on edges to $H_{1}$, then $|B| \leq 2t - 1$.
\end{itemize}
\end{claim}

\begin{proof}
For Item~$(1)$, suppose $|A| \geq t$. In order to avoid a monochromatic copy of $S_{t}^{2}$, there can be no two disjoint red edges or two disjoint blue edges within $A$ and no vertex within $A$ can have red or blue degree at least $2$. This means that if $A$ contains at least $2$ parts, each part must have order $1$ and there can be at most $2$ of them, a contradiction since $t \geq 3$. Thus, $A$ contains at most one part so $|A| \leq t - 1$, again a contradiction.

For Item~$(2)$, suppose $|B| \geq 2(t - 1) + 1$ so $B$ contains at least $3$ parts. Suppose all edges from $B$ to $H_{1} \cup H_{2}$ are red (meaning that all edges from $H_{1}$ to $H_{2}$ must be blue). In order to avoid a red copy of $S_{t}^{2}$, there can be at most one red edge within $B$. If $B$ has red edges between any pair of parts, those parts must have order $1$ each, leading to at least $4$ parts in $B$. In either case, all remaining edges between parts in $B$ must be blue, leading to a blue copy of $S_{t}^{2}$.

For Item~$(3)$, suppose $|B| \geq 2t - 1$, say with all red edges to $H_{1}$. To avoid a red copy of $S_{t}^{2}$, $B$ contains no red copy of $2K_{2}$. By deleting at most one vertex, the set $B$ satisfies all the conditions of Item~$(2)$ above so the proof is complete.
\end{proof}

The following fact is immediate and will be used later.

\begin{fact}\labelz{Fact:Fact4}
Let $G$ be a $3$-colored complete graph of order $t$ at least $6$, say using red, blue, and green. If the red subgraph of $G$ is contained within a triangle and the blue subgraph of $G$ is contained within a triangle and all other edges must be green, then $G$ contains a green copy of $S_{t}^{2}$.
\end{fact}

We break the proof into two main cases based on the value of $k$. 

\setcounter{case}{0}
\begin{case}
$k = 3$.
\end{case}

In this case, $n=5t$ so first suppose $2\leq m\leq 3$. By the minimality of $m$, we may assume $m=2$, say with corresponding parts $H_{1}$ and $H_{2}$. Suppose all edges between $H_{1}$ and $H_{2}$ are blue and without loss of generality, that $|H_{1}|\geq |H_{2}|$. Then $|H_{1}|\geq \lceil\frac{5t}{2}\rceil$. To avoid a blue copy of $S_{t}^{2}$, $H_{1}$ contains no two disjoint blue edges. By deleting at most two vertices from $H_{1}$, we can obtain a subgraph $H_{1}'$ with no blue edge. Clearly, $|H'_1|\geq \lceil\frac{5t}{2}\rceil-2>2t-1$ and $H_{1}'$ uses only two colors, so $H_{1}'$ contains a monochromatic $S_{t}^{2}$, a contradiction. We may therefore assume that $m\geq 4$.

Then we have the following claim.

\begin{claim}\labelz{Claim:Claim2}
$r\leq 2$.
\end{claim}

\begin{proof}
Assume, to the contrary, that $r\geq 3$. For any choice of $3$ large parts, say $H_{1}, H_{2}, H_{3}$, to avoid a monochromatic copy of $S_{t}^{2}$, the triangle in the reduced graph corresponding to the parts $H_{1}, H_{2}, H_{3}$ must not be monochromatic. Without loss of generality, we suppose the edges from $H_{2}$ to $H_{3}$ are red and all other edges between those parts are blue. To avoid a red copy of $S_t^2$, there is no vertex incident two red edges and no red
$2$-matching within $H_{1}$ or within $H_{2} \cup H_{3}$, and hence there is at most one red edge in $H_2$. Similarly, there is at most one blue edge within $H_2$. From Fact~\ref{Fact:Fact4}, if $t \geq 6$, then there is a green $S_t^2$ in $H_2$, a contradiction. Otherwise if $t \leq 5$, the result is still easy to verify.
\end{proof}

We now consider the following subcases based on the value of $r$.

\begin{subcase}
$r=2$.
\end{subcase}

Without loss of generality, suppose all edges from $H_{1}$ to $H_{2}$ are blue. To avoid a blue copy of $S_{t}^{+}$, there is no part with blue edges to both $H_{1}$ and $H_{2}$. Let $A$ be the set of parts with red edges to $H_{1}$ and blue edges to $H_{2}$, let $B$ be the set of parts with blue edges to $H_{1}$ and red edges to $H_{2}$, and let $C$ be the set of parts with red edges to $H_{1}\cup H_{2}$.

First suppose $C\neq \emptyset$. To avoid a blue copy of $S_t^2$, $H_{1} \cup H_{2}$ contains no blue $2K_{2}$ (avoiding the edges between the two parts) and there is no vertex incident two blue edges (within either part), and hence there is at most one blue edge within $H_{1}$ or $H_{2}$. Since $C \neq \emptyset$, there is no red $2K_2$ within $H_{1}\cup H_{2}$, and hence at least one of $H_1$ or $H_2$ contains no red edge. Without loss of generality, we assume that $H_1$ contains no red edge. With at most one edge within $H_{1}$ that is not green, there is a green $S_t^2$ within $H_{1}$, a contradiction.

Thus, we may assume $C=\emptyset$. From Claim~\ref{Claim:Claim1}, we have $|A|\leq t-1$, and $|B|\leq t-1$. To avoid a red copy of $S_{t}^{2}$, $H_{1}$ and $H_{2}$ each contain no red $2K_{2}$. Thus, by removing at most $4$ vertices from $H_{1} \cup H_{2}$, we can obtain a subgraph $H'$ with no red edges. Since $|H'| \geq |H_{1} \cup H_{2}| - 4 \geq 5t - 2(t - 1) - 4 = 3t - 2$, this subgraph contains either a blue or green copy of $S_{t}^{2}$, a contradiction.

\begin{subcase}
$r=1$.
\end{subcase}

Let $A$ be the set of parts with blue edges to $H_{1}$ and $B$ be the set of parts with red edges to $H_{1}$. If either of these sets is empty, then this contradicts the minimality of $m$ so $A, B \neq \emptyset$. To avoid a monochromatic copy of $S_{t}^{2}$, we see that $H_{1}$ contains no red or blue copy of $2K_2$. By Claim~\ref{Claim:Claim1}~$(2)$, we get $|A|, |B| \leq 2t - 2$. This means that $|H_{1}| \leq 5t - 2(2t - 2) = t + 4$. By removing at most $4$ vertices from $H_{1}$, we can obtain a subgraph $H'$ in which there are no red or blue edges. This subgraph having order at least $t$ and colored entirely in green means there is a green copy of $S_{t}^{2}$, a contradiction.



\begin{subcase}
$r=0$.
\end{subcase}

Let $A$ be the set parts with blue edges to $H_{1}$, and $B$ be the set parts with red edges to $H_{1}$. Note that by minimality of $t$, we have $A\neq \emptyset$ and $B\neq \emptyset$. Without loss of generality, suppose $|A|\geq |B|$. As $k = 3$, $n = 5t$, so $|A|\geq 2t+1$. To avoid a blue copy of $S_{t}^{2}$, $A$ contains no blue copy of $2K_{2}$. By deleting at most $2$ vertices from $A$, we can obtain a subgraph with no blue edges, call it $A'$. Then $|A'|\geq 2t-1$, and $A'$ contains at least $3$ parts since each part of $A'$ has order at most $t-1$. All edges in between the parts of $A'$ must be red (since there are no blue edges), but this yields a red copy of $S_{t}^{2}$, a contradiction.

\begin{case}
$k \geq 4$.
\end{case}

First suppose $2\leq m\leq3$. If $m\leq3$, then by the minimality of $m$, we may assume $m=2$, say with corresponding parts $H_{1}$ and $H_{2}$. Suppose all edges in between $H_{1}$ and $H_{2}$ are blue. To avoid a blue copy of $S_{t}^{2}$, the two parts $H_{1}$ and $H_{2}$ together contain no $2K_{2}$ in blue within the parts. Be deleting at most $2$ vertices from $H_{1}\cup H_{2}$, we can obtain a subgraph $H'$ which contains no blue edge. Appling induction on $k$ within the two parts of $H'$, this means that
$$
|G|=|H_{1}|+|H_{2}|\leq |H'| + 2 \leq 2[gr_{k-1}(K_{3}:S_{t}^{2})-1]+2<n,
$$
a contradiction.

If $r\geq 5$ and $m\geq 6$, then any choice of $6$ parts containing the $5$ parts $\mathcal{H}=\{H_{1}, H_{2}, H_{3}, H_{4},H_{5}\}$ will contain a monochromatic triangle in the reduced graph. Such a triangle must contain at least two parts from $\mathcal{H}$, meaning that the corresponding subgraph of $G$ must contain a monochromatic copy of $S_{t}^{2}$, a contradiction. Thus, we may assume either $4\leq m\leq5$ or $r\leq 4$. We consider the following subcases based on the value of $r$.

\begin{subcase}
$r=0$.
\end{subcase}

Let $A$ be the set of vertices with blue edges to $H_{1}$, and $B$ be the set of vertices with red edges to $H_{1}$. Note that by minimality of $t$, we have $A\neq \emptyset$ and $B\neq \emptyset$. Without loss of generality, suppose $|A|\geq |B|$. Since $k\geq 4$, we have $n\geq 10t$, so 
$$
|A|\geq \ceil{\frac{9t+1}{2}}-2 > 2t.
$$
On the other hand, by Claim~\ref{Claim:Claim1}, we have $|A| \leq 2t - 1$, a contradiction.


\begin{subcase}
$r=1$.
\end{subcase}

Let $A$ be the set of parts with blue edges to $H_{1}$ and $B$ be the set of parts with red edges to $H_{1}$. By Claim~\ref{Claim:Claim1}~$(2)$, we have $|A|\leq 2(t-1)$ and $|B|\leq 2(t-1)$. To avoid a monochromatic $S_{t}^{2}$, $H_{1}$ contains no blue or red $2K_{2}$. By deleting at most $4$ vertices from $H_{1}$, we can obtain a subgraph $H'$ with no blue or red edge. By applying induction on $k$ within $H'$, we get
$$
|G|=|H_{1}|+|A|+|B| \leq |H'| + 4 + 4(t - 1) \leq [gr_{k-2}(K_{3},S_{t}^{2})-1]+4t < n,
$$
a contradiction.

\begin{subcase}
$r=2$.
\end{subcase}

Suppose blue is the color of the edges between $H_{1}$ and $H_{2}$. There is no part with blue edges to $H_{1}\cup H_{2}$, since otherwise such a part would form a blue triangle in the reduced graph using the two large parts, a contradiction. Let $A$ be the set of vertices with blue edges to $H_{1}$ and red edges to $H_{2}$, $B$ be the set of vertices with red red edges to $H_{1}\cup H_{2}$, and $C$ be the set of vertices with red edges to $H_{1}$ and blue
edges to $H_{2}$. By Claim~\ref{Claim:Claim1}, we have $|A|\leq t-1$, $|B|\leq 2(t-1)$, and $|C|\leq t-1$. To avoid a blue copy of $S_{t}^{+}$, there can be at most one blue edge within $H_{1}$ or $H_{2}$. Since $H_{1}$ and $H_{2}$ must have red edges to some other part, each contains no red $2K_{2}$. By removing at most $2$ vertices from each part $H_{1}$ and $H_{2}$, we can remove all red edges from $H_{1} \cup H_{2}$ and by removing at most one vertex, we can remove all blue edges from within $H_{1}$ and $H_{2}$. Thus, $|H_{1} \cup H_{2}| \leq 2[gr_{k - 1}(K_{3}, S_{t}^{2}) - 1] + 5$, so we get 
$$
|G|=|H_{1}|+|H_{2}|+|A|+|B|+|C|\leq 2[gr_{k-2}(K_{3},S_{t}^{2})-1]+5+4(t-1)<n,
$$
a contradiction.

\begin{subcase}
$r=3$.
\end{subcase}

To avoid a monochromatic $S_{t}^{2}$, the triangle in the reduced graph corresponding to the parts $H_{1}, H_{2}, H_{3}$ must not be monochromatic. Without loss of generality, suppose the edges from $H_{2}$ to $H_{3}$ are red and all edges between $H_{1}$ and $H_{2} \cup H_{3}$ are blue. First we claim that there is no part with blue edges to $H_{1}$. Otherwise suppose there is one part, say $H_{4}$ with blue edges to $H_{1}$. Then to avoid a blue $S_{t}^{2}$, all edges from $H_{4}$ to $H_{2}\cup H_{3}$ are red, but then the red triangle  $H_{4}H_{2}H_{3}$ contains a red $S_{t}^{2}$, a contradiction. Thus all edges from $H_{1}$ to $H_{4}\cup \ldots \cup H_{m}$ are red and no vertex has red edges to $H_{2} \cup H_{3}$. Let $A$ be the set of vertices with blue edges to $H_{2}$ and with red edges to $H_{3}$, $B$ be the set of vertices with blue edges to $H_{2}\cup H_{3}$, $C$ be the set of parts with red edges to $H_{2}$ and with
blue edges to $H_{3}$. To avoid a monochromatic copy of $S_{t}^{2}$, $H_{1}\cup H_{2}\cup H_{3}$ contains at most one blue edge within a part, $H_{2}\cup H_{3}$ contains at most one red edge within the parts, and $H_{1}$ contains no red $2K_{2}$. Thus, by deleting at most $4$ vertices from $H_{1}\cup H_{2}\cup H_{3}$, we can obtain a subgraph $H'$ which contains no blue and red edge within the subgraphs induced on $H_{1}$, $H_{2}$, $H_{3}$. By Claim~\ref{Claim:Claim1}, $|A|\leq t-1$, $|B|\leq 2(t-1)$ and $|C|\leq t-1$. Then 
$$
|G|=|H_{1}|+|H_{2}|+|H_{3}|+|A|+|B|+|C|\leq 3[gr_{k-2}(K_{3},S_{t}^{2})-1]+4+4(t-1)<n,
$$
a contradiction.

\begin{subcase}
$r=4$.
\end{subcase}

Considering the subgraph of the reduced graph induced by the $r$ large parts, there can be no monochromatic triangle. This means there are two possible colorings. For the first coloring, suppose that all edges from $H_1$ to $H_2$ are blue and all edges from $H_3$ to $H_4$ are blue, and all other edges among these parts are red. For the second coloring, suppose that all edges from $H_1$ to $H_2\cup H_3$ are blue, all edges from $H_2$ to $H_4$ are blue, and all other edges among these parts are red.

For the first coloring, to avoid a monochromatic copy of $S_{t}^{2}$, we know that $t=4$, and the parts $H_{1}\cup H_{2}\cup H_{3}\cup H_{4}$ contain a total of at most one blue edge and one red edge within the parts. By deleting at most $2$ vertices from $H_{1}\cup H_{2}\cup H_{3}\cup H_{4}$, we can obtain a subgraph which contains no red or blue edge. Then 
$$
|G|=|H_{1}|+|H_{2}|+|H_{3}|+|H_{4}|\leq 4[gr_{k-2}(K_{3},S_{t}^{2})-1]+2<n,
$$
a contradiction.

For the second coloring, to avoid a monochromatic triangle, all edges from any part outside $\{H_1,H_2,H_3,H_4\}$ to $H_1\cup H_2$ must be red, and all edges from any part outside $\{H_1,H_2,H_3,H_4\}$ to $H_3\cup H_4$ must be blue. By minimality of $m$, this means there can be at most one such part, say $A$ and note that $|A| \leq t - 1$. As above, the parts $H_{1}\cup H_{2}\cup H_{3}\cup H_{4}$ contain a total of at most one blue edge and one red edge. By deleting at most $2$ vertices from $H_{1}\cup H_{2}\cup H_{3}\cup H_{4}$, we can obtain a subgraph which contains no red or blue edge. Then 
$$
|G|=|H_{1}|+|H_{2}|+|H_{3}|+|H_{4}|\leq 4[gr_{k-2}(K_{3},S_{t}^{2})-1]+2+t-1<n,
$$
a contradiction.

\begin{subcase}
$r=5$.
\end{subcase}

The reduced graph graph restricted to the $5$ large parts must be the unique $2$-coloring of $K_{5}$ containing two complementary copies of $C_{5}$. To avoid a monochromatic copy of $S_{t}^{2}$, $H_{1}\cup H_{2}\cup H_{3}\cup H_{4}\cup H_{5}$ can contain at most one blue edge within a part and one red edge within a part. By deleting at most $2$ vertices we can obtain a subgraph of $G$ such that none of the five parts contains a red or a blue edge. Then 
$$
|G|=|H_{1}|+|H_{2}|+|H_{3}|+|H_{4}|+|H_{5}|\leq 5[gr_{k-2}(K_{3},S_{t}^{2})-1]+2<n,
$$
a contradiction. \qed

\subsection{For general $r$ and $t$}

We first prove the following lemma which will be used later in this section. 

\begin{lemma}\labelz{Lem:Str1}
Let $G$ be a Gallai colored complete graph with parts
$H_1,H_2,\ldots,H_m$ such that $|H_i|\leq t-1$ for $1\leq i\leq m$.
If the edges in between each pair of parts $H_i,H_j \ (1\leq i\neq j\leq m)$ are
red and $\sum_{i=1}^m|H_i|\geq 2t+r$, then there is a red $S_t^r$ in
$G$.
\end{lemma}
\begin{proof}
Since $\sum_{i=1}^m|H_i|\geq 2t+r$, it follows that there are at
least $3$ parts in $G$. Without loss of generality, let $|H_1|\leq
|H_2|\leq \ldots \leq |H_m|$. Let $v\in H_1$. The number of red
edges from $v$ to $\sum_{i=2}^m|H_i|$ is at least $t+r+1$. Since
$|H_1|\leq |H_2|\leq t-1$, it follows that $\sum_{i=3}^m|H_i|\geq
r+2$, and hence there is a red $S_t^r$, a contradiction.
\end{proof}

The lower bound in Theorem~\ref{Thm:Str} is provided by Lemma~\ref{Lem:Lower-St2}, and the upper bound of Theorem~\ref{Thm:Str} is provided by the following.

\begin{lemma}\labelz{Lem:Str-Up}
For $k\geq 1$ and $t\geq 6r-5$,
$$
gr_k(K_3:S_t^{r})\leq \begin{cases}
[2t+8(r-1)]\times 5^{\frac{k-2}{2}}-4(r-1), &\mbox {\rm if}~k~is~even;\\[0.2cm]
[t+4(r-1)]\times 5^{\frac{k-1}{2}}-4(r-1), &\mbox {\rm if}~k~is~odd.
\end{cases}
$$
\end{lemma}

\begin{proof}
From Item $(3)$ of Theorem~\ref{Thm:St2Classic}, we have $R(S_t^{r},S_t^{r})=2t+2r-1\leq 2t+4r-4$, and hence the result is true for $k=2$. We therefore suppose $k\geq 3$ and let $G$ be a Gallai coloring of $K_n$ where
$$
n=n(k,r,t) =\begin{cases}
[2t+8(r-1)]\times 5^{\frac{k-2}{2}}-4(r-1), &\mbox {\rm if}~k~is~even;\\[0.2cm]
[t+4(r-1)]\times 5^{\frac{k-1}{2}}-4(r-1), &\mbox {\rm if}~k~is~odd.
\end{cases}
$$

Since $G$ is a Gallai coloring, it follows from Theorem~\ref{Thm:G-Part} that there is a Gallai partition of $V(G)$. Suppose that the two colors appearing in the Gallai partition are red and blue. Let $m$ be the number of parts in this partition and choose such a partition where $m$ is minimized. Let $H_{1}, H_{2}, \dots, H_{m}$ be the parts of this partition, say with $|H_{1}| \geq |H_{2}| \geq \dots \geq |H_{m}|$. When the context is clear, we also abuse notation and let $H_{i}$ denote the vertex of the reduced graph corresponding to the part $H_{i}$.

If $2\leq m\leq 3$, then by the minimality of $m$, we may assume $m=2$. Let $H_1$ and $H_2$ be the corresponding parts. Suppose all edges from $H_1$ to $H_2$ are red. If $|H_{i}|\geq t \ (i=1,2)$, then to avoid creating a red copy of $S_t^{r}$, there are at most $r-1$ disjoint red edges in each $H_i$ with $i=1,2$. If we remove all the vertices of a maximum set of disjoint red edges within $H_1$ and $H_2$, we create subgraphs $H_{1}'$ and $H_{2}'$ containing no red edge within either $H_{i}'$. This means that
$$
|G| = |H_{1}| + |H_{2}| \leq 2(n(k-1,r,t)-1)+2(2r-2)<n,
$$
a contradiction. If $|H_{1}|\geq t$ and $|H_{2}|\leq t-1$, then we similarly remove all the vertices of a maximum set of disjoint red edges within $H_1$ to create a subgraph $H_{1}'$ containing no red edge within $H_{1}'$. This means that
$$
|G| = |H_{1}| + |H_{2}| \leq n(k-1,r,t)-1+(2r-2)+t-1<n,
$$
a contradiction. Finally if $|H_{1}|\leq t-1$ and $|H_{2}|\leq t-1$, then $|G| = |H_{1}| + |H_{2}| \leq 2(t-1)<n$, a contradiction. Thus, we may assume $m \geq 4$ and by minimality of $m$, each part has edges to some other parts in both red and blue. If a part has order at least $t$, it can therefore contain no set of $r$ independent edges in either red or blue. By removing the at most $4r - 4$ vertices of a red maximum matching and a blue maximum matching from such parts, we can obtain a subgraph with no red or blue edges. This leads to the following fact.

\begin{fact}\labelz{Fact:NotTooBig}
$$
|H_{i}| \leq n(k - 2, r, t) - 1 + (4r - 4).
$$
\end{fact}

Let $a$ be the number of parts of the Gallai partition with order at least $t$ and call these parts ``large'' while other parts are called ``small''. Then $|H_a|\geq t$ and $|H_{a+ 1}|\leq t-1$. To avoid a monochromatic copy of $S_t^{r}$, there can be no monochromatic triangle within the reduced graph restricted to these $a$ large parts, leading to the following immediate fact.

\begin{fact}\labelz{Fact:n1}
$a\leq 5$.
\end{fact}

Let $A$ be the set of parts with red edges to $H_1$, and $B$ be the set of parts with blue edges to $H_1$. Then we have the following claim.

\begin{claim}\labelz{Claim16}
If $a \leq 1$, we have $|A|, |B| \leq 2t+3r-3$.
\end{claim}

\begin{proof}
Assume, to the contrary, that $|A|\geq 2t+3r - 2$. To avoid a red copy of $S_{t}^{r}$ centered in $H_{1}$, the subgraph $A$ must contain at most $r - 1$ independent red edges. By deleting all the vertices of a maximum red matching within $A$, we create a subgraph $A'$ containing no red edge with $|A'|\geq |A| - 2(r - 1) \geq 2t+r$. From Lemma~\ref{Lem:Str1}, $A$ contains a blue copy of $S_t^r$, a contradiction. The same holds for $B$.
\end{proof}

The remainder of the proof is broken into cases based on the value of $a$.

\setcounter{case}{0}
\begin{case}\labelz{Case:n1}
$a=0$.
\end{case}

Then by Claim~\ref{Claim16}, we have
$$
|G| \leq  |H_{1}'| + |A| + |B|  \leq (t-1) + 2(2t+3r-3) < n,
$$
a contradiction.

\begin{case}
$a=1$.
\end{case}

By Fact~\ref{Fact:NotTooBig} and Claim~\ref{Claim16}, we have 
$$
|G| = |H_{1}| + |A| + |B| \leq [n(k - 2, s,t) - 1 + 4(r - 1)] + 2[(2t+3r-3) + 4(r - 1)] < n,
$$
a contradiction.

\begin{case}\labelz{Case:nr5}
$a=5$.
\end{case}

In this case, $m=5$ since otherwise any monochromatic triangle in the reduced graph restricted to $H_{1}, H_{2}, \dots, H_{6}$ would yield a monochromatic copy of $S_t^r$. To avoid the same construction, the reduced graph on the parts $H_1,H_2,H_3,H_4,H_5$ must be the unique $2$-coloring of $K_{5}$ with no monochromatic triangle, say with $H_1H_2H_3H_4H_5H_1$ and $H_1H_3H_5H_2H_4H_1$ making two monochromatic cycles in red and blue respectively.

In order to avoid a red copy of $S_t^r$, it must be the case that the subgraph induced on $H_1\cup H_3$ contains at most $r-1$ disjoint red edges. Similarly $H_{1} \cup H_{4}$, $H_{2} \cup H_{4}$, $H_{2} \cup H_{5}$, and $H_{3}\cup H_{5}$ each contain at most $r- 1$ disjoint red edges. Putting these together, there are at most a total of $\frac{5r- 5}{2}$ disjoint red edges within the parts $H_{1}, H_{2}, \dots, H_{5}$. Thus, by deleting at most $5r-5$ vertices, we can obtain a subgraph of $G$ in which the $5$ parts contain no red edges and symmetrically, by deleting at most another $5r- 5$ vertices, we can obtain a subgraph $G' \subseteq G$  in which the $5$ parts also contain no blue edges. This means that
$$
|G| \leq |G'| + (10r-10) \leq 5[n(k - 2, r,t) - 1] + (10r - 10) < n,
$$
a contradiction.

\begin{case}\labelz{Case:r4}
$a=4$.
\end{case}

To avoid monochromatic triangle within the reduced graph restricted to the four large parts, these parts must form a red path, say $H_1H_2H_3H_4$, and a blue path, say $H_2H_4H_1H_3$. Since there are at most $r-1$ independent red edges within $H_{1} \cup H_{3}$, by deleting at most $2r-2$ vertices from $H_1\cup H_3$, we can obtain a subgraph of $H_{1} \cup H_{3}$ in which there are no red edges. Similarly, by deleting at most $2r-2$ vertices from $H_2\cup H_4$, we can obtain a subgraph of $H_{2} \cup H_{4}$ in which there are no red edges. Symmetrically, if we delete at most $4r-4$ vertices in $H_1\cup H_2\cup H_3\cup H_4$, there are no blue edges within $H_{1} \cup H_{2}$ or within $H_{3} \cup H_{4}$. If $m = 4$, this means that
$$
|G|\leq 4[n(k - 2, r,t) - 1]+ (8r-8) <n,
$$
a contradiction. If $m > 4$, then to avoid a monochromatic triangle in the reduced graph that includes at least two large parts, all small parts must have red edges to $H_{1} \cup H_{4}$ and blue edges to $H_{2} \cup H_{3}$. Then by minimality of $m$, we must have $m \leq 5$. Now we may apply the same argument as in Case~\ref{Case:nr5} to complete the proof in this case.

\begin{case}
$a=3$.
\end{case}

The triangle in the reduced graph corresponding to the three large parts cannot be monochromatic so without loss of generality, suppose the edges from $H_{1}$ to $H_{2} \cup H_{3}$ are red, and $H_2H_3$ is blue. To avoid a red or blue triangle, any remaining parts are partitioned into the following three sets.
\begin{itemize}
\item Let $A$ be the set of parts outside $H_1,H_2,H_3$ each with all blue edges to $H_1,H_3$ and all red edges to $H_2$,
\item let $B$ be the set of parts outside $H_1,H_2,H_3$ each with all red edges to $H_2,H_3$ and all blue edges to $H_1$, and
\item let $C$ be the set of parts outside $H_1,H_2,H_3$ each with all blue edges to $H_1,H_2$ and all red edges to $H_3$.
\end{itemize}

Note that $A \cup B \cup C \neq \emptyset$ (recall each part must have red and blue edges to some other parts). By Fact~\ref{Fact:NotTooBig}, we have $|H_1|+|H_2|+|H_3|\leq 3[n(k - 2, r,t) - 1]+12r-12$. Since $H_1A,H_1B,H_1C$ are blue, it follows that by deleting at most $2r-2$ vertices from $A\cup B\cup C$, we can obtain a subgraph in which there are no blue edges. Let $A'\cup B'\cup C'$ be this remaining graph. Then the edges in between the parts within $A'\cup B'\cup C'$ are all red. Since $H_2A,H_2B$ are blue, it follows that by deleting at most $2r-2$ vertices in $A'\cup B'$, we can obtain a subgraph in which there are no red edges and so there is only one (small) part. This means that $|A'|+|B'|\leq (t-1)+(2r-2)$. By deleting at most $2r-2$ vertices from $C'$, we can obtain a subgraph in which there are no red edges and so there is only one (small) part. This means $|C'|\leq (t-1)+(2r-2)$, implying that $|A'|+|B'|+|C'|\leq 2(t-1)+(4r-4)$, and therefore $|A|+|B|+|C|\leq 2(t-1)+(6r-6)$. This gives us
\beqs
|G| & = & |A|+|B|+|C|+|H_1|+|H_2|+|H_3|\\
~ & \leq & 2(t - 1) + (6r - 6) + 3[n(k - 2, r, t) - 1] + (12r - 12)\\
~ & < & n, 
\eeqs
a contradiction.

\begin{case}
$a=2$.
\end{case}

To avoid creating a monochromatic copy of $S_t^r$, there is no part outside $H_1$ and $H_2$ with red edges to all of $H_1\cup H_2$. Let $A$ be the set of vertices outside $H_1 \cup H_2$ each with all blue edges to $H_2$ and all red edges to $H_1$, and let $B$ be the set of vertices outside $H_1 \cup H_2$ each with all blue edges to $H_1 \cup H_2$, and let $C$ be the set of vertices outside $H_1 \cup H_2$ each with all blue edges to $H_1$ and all red edges to $H_2$. As in the previous case, we get $|A|\leq (t-1)+(4r-4)$ and $|C|\leq (t-1)+(4r-4)$. If $|B|\geq t$, then there are at most $r-1$ disjoint blue edges within $B\cup C$ or $B\cup A$. By deleting at most $2r-2$ vertices in $B\cup C$ or $B\cup A$, we can obtain a subgraph in which there are no blue edges. By Lemma~\ref{Lem:Str1}, this subgraph has order at most $2t+r-1$, and so $|B\cup C|\leq 2t+3r-3$ and $|B\cup A|\leq 2t+3r-3$. Therefore, again using Fact~\ref{Fact:NotTooBig}, we have 
\beqs
|G| & = & |A|+|B|+|C|+|H_1|+|H_2|\\
~ & \leq & (t-1)+(4r-4)+(2t+3r-3)+2[n(k - 2, s,t) - 1]+(8r-8)\\
~ & < & n,
\eeqs
a contradiction. If $|B|\leq t-1$, then 
\beqs
|G| & = & |A|+|B|+|C|+|H_1|+|H_2|\\
~ & \leq & 3(t-1)+2(4r-4)+2[n(k - 2, s,t) - 1]+(8r-8)\\
~ & < & n, 
\eeqs 
a contradiction.
\end{proof}

\end{document}